\documentclass[10pt,twocolumn,twoside] {IEEEtran}

\usepackage{times} 
\usepackage{pdfsync} 
\usepackage{epsfig}
\usepackage{epstopdf}
\usepackage[cmex10]{amsmath}
\usepackage{amsmath} 
\usepackage{amssymb}
\usepackage{enumerate}
\usepackage{color} 
\usepackage{cite} 
\usepackage{array}
\usepackage{mathrsfs}
\usepackage{fix2col} 
\usepackage{hyperref}

\graphicspath{{TCNS-13-0019-Fig/}}

\newcommand{\Prob}{\mathbb{P}}
\newcommand{\E}{\mathbb{E}}

\def\qed
{\hfill\vbox{\hrule width 0.5em\nointerlineskip\hbox to
0.5em{\vrule height 0.5em \hfill\vrule height
0.5em}\nointerlineskip\hrule width 0.5em}}

\newtheorem{thm}{Theorem}[section]
\newtheorem{defi}{Definition}[section]

\newtheorem{lem}{Lemma} [section]

\newtheorem{remark}{Remark}[section] 

\title{An Optimal Transmission Strategy for Kalman Filtering over Packet Dropping Links with Imperfect Acknowledgements}

\author{Mojtaba Nourian, Alex S. Leong, Subhrakanti Dey and Daniel E. Quevedo
\thanks{M. Nourian and A. S. Leong are with the Department of Electrical and Electronic Engineering, The University of Melbourne, VIC 3010, Australia (email: {\tt\small \{mojtaba.nourian,asleong\}@unimelb.edu.au}). 

S. Dey is with the Department of Engineering Sciences, Uppsala University, Sweden (email: {\tt\small subhra.dey@signal.uu.se}).

D. E. Quevedo is with the School of Electrical Engineering \& Computer Science, The University of Newcastle, NSW 2308, Australia (email: {\tt\small dquevedo@ieee.org}).}
}

\begin{document}

\maketitle

\begin{abstract} This paper presents a novel design methodology for optimal transmission policies at a smart sensor to remotely estimate the state of a stable linear stochastic dynamical system. The sensor makes measurements of the process and forms estimates of the state using a local Kalman filter. The sensor transmits  quantized information over a packet dropping link to the remote receiver. The receiver sends packet receipt acknowledgments back to the sensor via an erroneous feedback communication channel which is itself packet dropping. The key novelty of this formulation is that the smart sensor decides, at each discrete time instant, whether to transmit a  quantized version of either its local state estimate or its local innovation. 
The objective is to design optimal transmission policies in order to minimize a long term average  cost function as a convex combination of the receiver's expected estimation error covariance and the energy needed to transmit the packets. The optimal transmission policy is obtained by the use of dynamic programming techniques. Using the concept of submodularity, the optimality of a threshold policy in the case of scalar systems with perfect packet receipt acknowledgments is proved. Suboptimal solutions and their structural results are also discussed. Numerical results are presented illustrating the performance of the optimal and suboptimal transmission policies.

\end{abstract}

\begin{IEEEkeywords} Wireless sensor networks, state estimation, packet drops, high resolution quantizer, Markov decision processes with imperfect state information, threshold policy.
\end{IEEEkeywords}

\section{Introduction} 


\IEEEPARstart{O}{ne} of the important challenges in wireless based networks is to improve system performance and reliability under resource (e.g., energy/power, computation and communication) constraints. This concern is particularly crucial in industrial applications such as remote sensing and real-time control where a high level of reliability is usually required. As a consequence, it becomes of significant importance to investigate the impact of realistic wireless communication channel models in the area of state estimation and control of networked systems \cite{baillieul2007control}. Two important limitations of wireless communication channels in these problem formulations include: (i) limited bandwidth, and (ii) information loss.

Among the many papers in the area of networked state estimation and control over bandwidth limited  channels, we first mention  \cite{nair2004stabilizability}, which addresses the minimum data rate required for stability of a linear stochastic system with quantized measurements received through a finite rate channel. Recently, this work is extended to the general case of time-varying Markov digital communication channels in \cite{Minero_13}. The reader is also referred to the survey \cite{nair2007feedback} and the references therein.

Since the seminal work of \cite{Sinopoli}, state estimation or Kalman filtering problems over packet dropping communication channels have been extensively studied (see for example\cite{XuHespanha,HuangDey,Epstein_Automatica,Schenato,MoSinopoli,Quevedo_Automatica}, among others). The reader is also referred to the survey \cite{Schenato_proceedings} and the references therein. In these problems sensor measurements (or state estimates in the case of \cite{XuHespanha}) are grouped into packets which are transmitted over a packet dropping link. The focus in these works is on deriving conditions on the packet arrival rate in order to guarantee the stability of the Kalman filter. There are other works which are concerned with estimation performance (e.g. minimizing the expected estimation error covariance) rather than just stability. For instance, power allocation techniques have been applied to the Kalman filtering problem in \cite{QuevedoAhlenOstergaard,Alex_CDC12,Nourian_TAC13} in order to improve the estimation performance and reliability.

Even though most of the works available in the literature focus on only one of the two mentioned communication limitations (limited bandwidth or information loss), some recent works attempt to address both limitations. In particular, the problem of minimum data rates for achieving bounded average state estimation error in linear systems over lossy channels is studied in \cite{you2010minimum,YouXie11} 
(see also 
\cite{Ishii13}), while the problem of state control around a target state trajectory in the case of both signal quantization and packet drops is investigated in \cite{Trivellato_10,QuevedoOstergaardNesic}. The work in \cite{DeyChiuso:CDC13} concentrates on designing coding and decoding schemes to remotely estimate the state of a scalar stable stochastic linear system over a communication channel subject to both quantization noise and packet loss. 

\begin{figure*}[!t]
\centering 
\includegraphics[height=3.5cm,width=14cm]{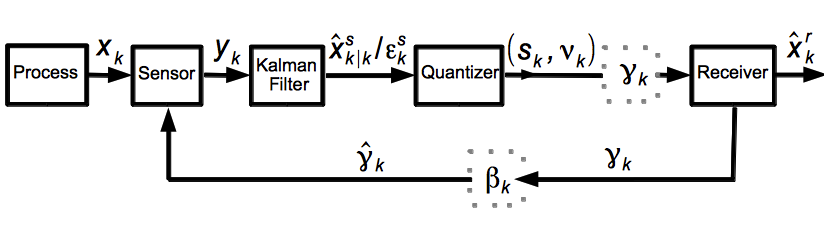} 
\caption{System Architecture}
\label{system_model}
\end{figure*}

Similar to \cite{DeyChiuso:CDC13}, the current paper is concerned with remote state estimation   subject to both quantization noise and packet drops. However, rather than considering fixed coding and decoding schemes, we are interested in choosing optimal transmission policies at the smart sensor that decides between sending the sensor's local state estimates or its local innovations. 
More specifically, we present a novel design methodology for optimal transmission policies at a smart sensor to remotely estimate the state of a stable (see Section \ref{conc_sec} for some comments on extensions to unstable systems) linear stochastic dynamical system. The sensor makes measurements of the process and forms estimates of the state using a local Kalman filter (see Fig. \ref{system_model}). The sensor then transmits quantized (using a high resolution quantizer) information over a packet dropping link to the remote receiver. The sensor  decides, at each  time instant, whether to transmit a  quantized version of either its local state estimate or its local innovation.  The receiver runs a Kalman filter with random packet dropouts to minimize the estimation error covariances based on received measurements. 

The packet reception probability is generally a function of the length of the packet, such that shorter packets (and hence lower required data rates) may result in higher packet receipt probabilities. Since the local innovation process has a smaller covariance, for a fixed packet reception probability the quantized innovations require less energy to transmit than the quantized state estimates. However, due to the packet dropping link between the sensor and the remote estimator, if there has been a number of successive packet losses then receiving a quantized state estimate might be more beneficial in reducing the estimation error covariance at the remote estimator than receiving the innovations. Thus, there is a tradeoff between whether the sensor should transmit its local state estimates or its local innovations.
In general, knowledge at the sensor of whether its transmissions have been received is achieved via some feedback mechanism. Here, in addition to the case of perfect packet receipt acknowledgments, we consider the more difficult problem where the feedback channel from the receiver to the sensor is an erroneous packet dropping link. 

The objective is to design optimal transmission policies in order to minimize a long term average (infinite-time horizon) cost function as a convex combination of the receiver's expected estimation error variance and the energy needed to transmit the packets. This problem is formulated as an average cost Markov decision process with imperfect state information. The optimal transmission policy is obtained by the use of dynamic programming techniques. 

In summary, the main contributions of this paper are as follows:

\begin{enumerate}[i)]

\item Unlike a large number of papers focusing on only one of the two communication limitations (limited bandwidth or information loss),
we consider  both limitations, i.e.,  remote state estimation   subject to both quantization noise and packet drops. 

\item Although recent work such as \cite{YouXie11, Ishii13} consider packet loss and data rate constraints simultaneously, the focus of these papers is on stabilizability (implying only bounded estimation error) whereas the focus on our work is on the  actual estimation error performance of the remote estimator (albeit for a stable system) and the optimization of a cost combining the long term average of estimation error and transmission energy expenditure. 

\item Unlike \cite{DeyChiuso:CDC13} which considers fixed coding and decoding schemes, we are interested in choosing optimal transmission policies at the smart sensor that decides between sending the sensor's local state estimates or its local innovations. 

\item We consider the case of imperfect feedback acknowledgements, which is more difficult to analyze than the case of perfect feedback acknowledgements. We model the feedback channel by a general erasure channel with errors.

\item It is well known that the optimal solution obtained by a stationary control policy minimizing the infinite horizon control cost is computationally prohibitive. Thus motivated, for the scalar case we provide structural results on the optimal policy which lead to simple threshold policies which are optimal and yet very simple to implement.
 
\item Finally, also motivated by the computational burden for the optimal control solution in the general case of imperfect acknowledgments, we provide a sub-optimal solution based on an estimate of the error covariance at the receiver. Numerical results are presented to illustrate the performance gaps between the optimal and sub-optimal solutions.

\end{enumerate}

The organization of the paper is as follows. The system model is given in Section \ref{sec:sensor}. The augmented state space model at the remote receiver is constructed in Section \ref{sec:aug} and the corresponding Kalman filtering equations are given. Section \ref{sec:opt} presents optimal transmission policy problems, together with their solutions, in both the cases of perfect and imperfect packet receipt acknowledgements. A suboptimal transmission scheme in the case of imperfect packet receipt acknowledgements is considered in Section \ref{sec:subopt}. For scalar systems, Section \ref{sec:thresh} proves  the optimality of the threshold transmission policy for the case of perfect packet receipt acknowledgements. Numerical simulations are given in Section \ref{numerical_sec}. 

\section{System Model} \label{sec:sensor}


We use the following notation. Let $(\Omega, \mathcal F,\Prob)$ be a complete probability space. $\E$ denotes the expectation. Throughout the paper, the subscript or superscript $s$ are used for the sensor's quantities, and the superscript $r$ is used for the receiver's quantities. We say that a matrix $X>0$ if $X$ is positive definite, and $X \geq 0$ if $X$ is positive semi-definite. 

A diagram of the system architecture is shown in Fig. \ref{system_model}. Detailed descriptions of each part of the system is given below.
\subsection{Process Dynamics and Sensor Measurements} \label{ProcessDyn_sec}

We consider a stable uncontrolled linear time-invariant stochastic dynamical process 
\begin{align}   \label{ProcDyn}
& x_{k+1} = A x_k+w_k, \qquad k \geq 0
\end{align}
where $x_k \in \mathbb{R}^n$ is the process state at instant $k \geq 0$, with $A$ being a Schur stable matrix, and $\{w_k : k \geq 0\}$ is a sequence of independent and identically distributed (i.i.d.) Gaussian noises with zero mean and covariance $\Sigma_w \geq 0$. The initial state of the process $x_0$ is a Gaussian random vector, independent of the process noise sequence $\{w_k: k \geq 0\}$, with mean $\bar x_0 := \mathbb{E}[x_0]$ and covariance $P_{x_0} \geq 0$. \\
\indent

The sensor measurements are obtained in the form
\begin{align} \label{ProcObs}
& y_k= C x_k+v_k, \qquad k \geq 0
\end{align}
where $y_k \in \mathbb{R}^m$ is the vector observation at instant $k \geq 0$, $C \in \mathbb{R}^{m \times n}$, and $\{v_k: k \geq 0\}$ is a sequence of i.i.d. Gaussian noises, independent of both  $x_0$ and  $\{w_k: k \geq 0\}$, with zero mean and covariance $\Sigma_v > 0$. We enunciate the following assumption: \\
({\bf A1}) We assume that $(A,C)$ is detectable.  \qed

\subsection{Local Kalman Filter at the Smart Sensor} \label{sec:wireless-sensor-Kalman}
We assume that the sensor has some computational capabilities. In particular, it can run a local Kalman filter to reduce the effects of measurement noise, as in e.g. \cite{XuHespanha}.

Denote the local sensor information at each instant $k$ by $\mathcal Y_k^s := \sigma\{y_t: 0 \leq t \leq k\}$, which is the $\sigma$-field generated by the sensor measurements up to time $k$. We use the convention $\mathcal Y_0^s :=\{\emptyset,\Omega\}$. Then, the Kalman filtering and prediction estimates of the process state $x_k$ at the sensor are given by $\hat{x}_{k|k}^s =  \mathbb{E} [x_{k} | \mathcal Y_k^s]$ and $\hat{x}_{k+1|k}^s =  \mathbb{E} [x_{k+1} |\mathcal Y_k^s]$, respectively. 

We assume that the local Kalman filter has reached steady-state. The stationary error-covariance is defined by $P_s = \lim_{k \rightarrow \infty} \E[(x_{k+1}-\hat{x}_{k+1|k}^s)(x_{k+1}-\hat{x}_{k+1|k}^s)^T|\mathcal Y_k^s]$, which is the solution of the algebraic Riccati equation (see e.g. \cite{andersonMoore})
\begin{align} \label{DARE}
& \!\! P_s = A P_s A^T+ \Sigma_w \!\!-\!\! A P_s C^T (CP_sC^T + \Sigma_v)^{-1} C P_s A^T.
\end{align}
The Kalman filter equations for $\hat{x}_{k|k}^s$ and $\hat{x}_{k+1|k}^s$ are given by
\begin{align}
&  \hat x^s_{k|k} = \hat x^s_{k|k-1} + K_f (y_k-C \hat x^s_{k|k-1}), \quad k \geq 0 \label{KFR:Filter}\\
& \hat x^s_{k+1|k} = A \hat x^s_{k|k-1} + K_s (y_k-C \hat x^s_{k|k-1}), \quad k \geq 0 \label{KFR:Pred}
\end{align}
with $\hat x^s_{0|-1}:= \bar x_0$, where $K_f:= P_sC^T(CP_sC^T+\Sigma_v)^{-1} $ and $K_s:=A K_f$ are the stationary Kalman filtering and prediction gains, respectively. Denote the covariance of the local state estimate via $\Sigma_s := \lim_{k \rightarrow \infty} \E[(\hat{x}_{k+1|k}^s)(\hat{x}_{k+1|k}^s)^T|\mathcal Y_k^s]$, which satisfies the stationary Lyapunov equation
\begin{align}   \label{DLE}
&\Sigma_s = A \Sigma_s A^T + K_s (CP_sC^T+ \Sigma_v) K_s^T. 
\end{align}

\subsection{Coding Alternatives at the Smart Sensor} \label{SSCS}

We define the innovation process\footnote{Note that $\epsilon^s_{k}$ is a linear transformed version of the true innovation process of Kalman filtering given by $y_k-C\hat x^s_{k|k-1}$.} at the sensor $\epsilon_{(\cdot)}$ as
\begin{align} \label{InnDef}
 & \epsilon_k^s = \hat x_{k|k}^s - \hat x_{k|k-1}^s = K_f (y_k-C \hat x_{k|k-1}^s), \quad k\geq 0.
\end{align}
As depicted in Fig. \ref{system_model}, the sensor communicates over a digital erasure channel with a remote receiver which utilizes the received data to calculate an estimate of the process state $x_{(\cdot)}$. 

This work aims to investigate what data the smart wireless sensor should transmit to the receiver. Motivated by differential Pulse-Code Modulation (PCM) techniques \cite{queahl09b,jayant1985digital}, the digital sensor may convey either a vector quantized version of its local estimate or a vector quantized version of its innovation. Therefore, we may denote the packet sent by the sensor as
\begin{align}
  \label{sensor_packet}
 & s_k :=
  \begin{cases}
    \hat x_{k|k}^s+q_k^x & \text{if $\nu_k=1$}\\
    \epsilon_k^s  + q_k^\epsilon &\text{if $\nu_k=0$}
  \end{cases}, \quad k \geq 0
\end{align}
where $\nu_k\in\{0,1\}$ is a decision variable which is transmitted to the receiver in addition to  $s_k$. The sequence $\{\nu_k\}$ is designed at the sensor, see Section \ref{sec:opt}. In (\ref{sensor_packet}), $q_{(\cdot)}^x$ and $q_{(\cdot)}^\epsilon$ are the quantization noises resulting from encoding $\hat x_{k|k}^s$ and $\epsilon_k^s$ respectively. We note that in this paper the effects of the quantizer are only modelled via the additive quantization noise
 term in (\ref{sensor_packet}), which is assumed to be zero-mean white noise processes independent of the 
quantized signal. For high-rate quantization, such an approach is quite accurate (see Remark \ref{lowratequant} below for the validity of this model to low-moderate rate quantization), since the quantization noises at high rates are approximately uncorrelated with the quantizer inputs \cite{marco2005validity,Goyal}. It is also reasonable to assume that the quantization noises, whilst uncorrelated to the inputs, have covariances which are proportional to the input covariances, i.e.,
\begin{equation}
\label{Quant-var}
\begin{split}
& \!\! \Sigma_q^x := \lim_{k\to\infty}\E [q_k^x (q_k^x)^T] =\alpha_1 \lim_{k\to\infty}\E [\hat x_{k|k}^s(\hat x_{k|k}^s)^T]\\
& \!\! \Sigma_q^\epsilon := \lim_{k\to\infty}\E [q_k^\epsilon (q_k^\epsilon)^T]  = \alpha_0 \lim_{k\to\infty} \E [\epsilon_k^s(\epsilon_k^s)^T]
\end{split}
\end{equation}
for given $\alpha_0, \alpha_1 \geq 0$ which depends upon the quantizers and the bit-rates used\footnote{For an explanation on how the scaling factors $\alpha_0, \alpha_1 \geq 0$ arise, see 
page 3860 of \cite{LeongDeyNair_kalmanquant_journal}.}. We can obtain $ \lim_{k\to\infty} \E[\hat x_{k|k}^s(\hat x_{k|k}^s)^T] = \Sigma_s + K_f (CP_sC^T+\Sigma_v) K_f^T$ from (\ref{KFR:Filter}), and $ \lim_{k\to\infty}\E[\epsilon_k^s (\epsilon_k^s)^T] = K_f (CP_sC^T+\Sigma_v) K_f^T$ from (\ref{InnDef}).  

Consider a vector Gaussian source $s$ with $N=2^n$ quantizer levels where $n$ is the transmission rate (i.e., the number of bits transmitted per sample).  Then the quantization noise covariance of a high resolution quantizer will be $\Sigma_q \approx \alpha \E [s s^T]$.  For the case of asymptotically optimal lattice vector quantizers with Voronoi cell $S_0$, we have (see \cite{Moo_thesis}) 
\begin{equation*}
\begin{split}
\alpha= \frac{M(S_0) V^{2/m}}{\eta^2}  \frac{\frac{2}{m} \ln N}{N^{2/m}}
\end{split}
\end{equation*}
where  $m$ represents the dimension of the vector to be quantized, $\eta =  \sqrt{1/2}$,
$V =  \frac{\pi^{m/2}}{\Gamma(m/2+1)}$, 
$$M(S_0) = \frac{\frac{1}{k} \int_{S_0} || x-y||_2^2 dx}{v(S_0)^{1+2/m}}$$ is the normalized moment of inertia of $S_0$, and $v(S_0)$  the volume of $S_0$. 
For $m=1$, it can be shown that $\alpha$ reduces to $\alpha = \frac{4 \ln N}{3N^2}$. 
For the case of ``optimal'' Lloyd-Max quantizers, we have 
$\alpha \sim \frac{B_m}{N^{2/m}}$ (see \cite{Gersho}). However, the exact values of the constants $B_m$ are not known for dimensions $m \geq 3$. 
For $m=1$, we have $\alpha = \frac{\pi \sqrt{3}}{2N^2}$. 
\begin{remark}
In principle, this additive white noise model for the quantization error is theoretically valid for high resolution quantization. However, it has been reported by many works including the seminal review paper by Gray and Neuhoff  \cite{GrayNeuhoff} (see p. 2358) that the high resolution theory is fairly accurate for rates greater than or equal to $3$ bits per sample per signal dimension. More recent papers such as \cite{LeongDeyNair_kalmanquant_journal} have reported similar results in designing decentralized linear estimation schemes with quantized innovations. Finally, the same quantization noise model has been used in a parallel work by Dey, Chiuso and Schenato (see the extended online version of \cite{DeyChiuso:CDC13}).  It has been shown in 
\cite{DeyChiuso:CDC13} that only 3  bits of quantization per sample for a convex combination of the (scalar) state estimate and the innovation signal at the transmitter achieve a remote estimation error performance that is sufficiently close to the one predicted by the additive white noise model.   Note that
in the context of modern wireless LANs, communication rates of the order of Mega  bits per second are quite common implying that 3-5 bits of quantization per sample 
can be easily achieved.  Thus this approximation is a fairly accurate tool for analysis that is suitable for practical implementations as well. 
\label{lowratequant}
\end{remark}
\begin{remark}
Although quantization noise is generally modelled as uniformly distributed, it has been also shown in a number of works that a Gaussian approximation to the quantization noise is  valid at high rate quantization. In particular, quantization noise due to lattice vector quantization (as used in this work) approaches a white Gaussian noise in a divergence sense \cite{zamirfeder1996} as the resolution increases. 
\label{gaussquant}
\end{remark}
Based on the above discussion, we model the quantization noise processes $q_k^x$ and $q_k^{\epsilon}$ as 
zero-mean additive white Gaussian noise processes with covariances $\Sigma_q^x,  \Sigma_q^{\epsilon}$ 
respectively. While this model is valid in principle at high rate quantization, it serves as a good approximation and a very useful analytical tool also at low to moderate rates of quantization as explained above.

In what follows,
we allow the sensor to choose a varying rate of quantization in order to make the traces of the quantization noise covariances $\Sigma_q^x$ and $\Sigma_q^\epsilon$ the same.
From (\ref{Quant-var}), this implies that the data rates $n_0$ and $n_1$ for  transmitting $\epsilon_k^s$ and $\hat x_{k|k}^s$ in the case of the lattice vector quantizer satisfy
\begin{align*}
& \textrm{Tr} \Sigma_q^x \equiv \frac{M(S_0) V^{2/m}}{\eta^2}  \frac{2n_1\ln2 /m  }{2^{2n_1/m}} \\
& \hspace{1cm} \times \textrm{Tr}\big(\Sigma_s + K_f (CP_sC^T + \Sigma_v) K_f^T\big ) = \textrm{Tr} \Sigma_q^\epsilon \\
& \hspace{0.7cm} \equiv \frac{M(S_0) V^{2/m}}{\eta^2} \frac{2n_0\ln2 /m }{2^{2n_0/m}} \textrm{Tr}\big(K_f (CP_sC^T + \Sigma_v) K_f^T\big)
\end{align*} 
and  in the case of the Lloyd-Max quantizer
\begin{align}
& \textrm{Tr} \Sigma_q^x \equiv \frac{B_m}{2^{2n_1/m}}  \textrm{Tr}\big(\Sigma_s + K_f (CP_s C^T+ \Sigma_v) K_f^T\big) \notag \\
& \hspace{0.75cm} = \textrm{Tr} \Sigma_q^\epsilon \equiv \frac{B_m}{2^{2n_0/m}}  \textrm{Tr}\big(K_f (CP_sC^T + \Sigma_v) K_f^T\big) . \label{LMQuzn}
\end{align} 
If the resulting $n_0$ and $n_1$ are not integers, their nearest integers will be chosen as the transmission rates.   Since $\Sigma_s \geq 0$, we have $n_0 \leq n_1$ in the two cases above. 

Since the local innovation process has a smaller stationary covariance, and hence a smaller data rate to maintain a given packet receipt probability, transmitting $\epsilon_k^s$ should require less energy than transmitting $\hat x_{k|k}^s$ (see Section \ref{sec:comm-chann}). However, due to the packet dropping link between the sensor and the remote estimator, if there has been a number of successive packet losses then receiving $\hat x_{k|k}^s$ might be more beneficial in reducing the estimation error covariance at the remote estimator than receiving $\epsilon_k^s$. 
Thus, in this model it is not immediately clear whether the sensor should transmit local estimates $\hat x_{k|k}^s$ or innovations $\epsilon_k^s$. The present work seeks to elucidate this dilemma in answering how to optimally design the control sequence $\{\nu_k: k \geq 0\}$ using causal information available at the sensor. 

\subsection{Forward Erasure Communication Channel} \label{sec:comm-chann}

We assume that the forward communication channel between the sensor and the receiver is unreliable, see Figure \ref{system_model}. This channel carries $\{(s_k,\nu_k): k \geq 0\}$ and is characterized by the transmission success process $\{\gamma_k: k \geq 0\}$, where $\gamma_k=1$ refers to successful reception of $(s_k,\nu_k)$ and $\gamma_k=0$ quantifies a dropout. Since the decision variable $\nu_k$ consists of only one bit of information, it can be easily sent along with $s_k$ as a header in the transmitted packet.

In this work we assume that $\gamma_k$ is a Bernoulli random variable with $\mathbb{P}(\gamma_k=1)= 1-p$, where $p\in [0,1]$ is the packet loss probability.  
The packet loss probability is generally a function of the data rates,  such that higher data rates result in higher packet loss probabilities. If $p_b$ is the error probability of sending one bit, then the packet loss probability of sending a packet of $n$ bits will be of the form 
\begin{align} \label{PLP}
& p = 1-(1-p_b)^n
\end{align}
where the packet is assumed to be lost if an error occurs in any of its bits (e.g. when there is no channel coding used).
We assume that the bit error probability $p_b$ of a wireless communication channel depends on the transmission energy per bit $E_b$ such that  $p_b$ decreases as $E_b$ increases. The bit error probability $p_b$ can be computed  for different combinations of channels and digital modulation schemes. For example, in the case of Additive White Gaussian Noise (AWGN) channel with Binary Phase-Shift Keying (BPSK) modulation:
\begin{align} \label{BEP}
& p_b = Q \big(\sqrt{\frac{2E_b}{N_0}}\big) 
\end{align}
where  $N_0/2$ is the noise power spectral density, and $Q(x):= (1/\sqrt{2\pi})\int_x^\infty e^{-t^2/2}dt=\frac{1}{2}\textrm{erfc}(\frac{x}{\sqrt 2}) $ is the $Q$-function \cite{Proakis}. As a consequence of (\ref{PLP}) and (\ref{BEP}), to obtain a fixed packet dropout probability, when innovations are sent the transmit energy per bit will be lower than when local estimates are transmitted. In Section \ref{sec:opt} we will further elucidate the situation and allocate power levels accordingly. 

\subsection{Erroneous Feedback Communication Channel} \label{ternary_subsec}

In the present work we will study the more realistic but complex case where
acknowledgments are unreliable (see \cite{garsin10,moayedi} for  relevant models with imperfect feedback mechanism). In this case, the packet loss process $\{\gamma_k, k \geq 0\}$ is not known to the sensor. Instead, the sensor receives an imperfect acknowledgment process $\{\hat \gamma_k, k \geq 0\}$ from the receiver. It is assumed that after the transmission of $y_k$ and before transmitting $y_{k+1}$ the sensor has access to the ternary process $\hat \gamma_k\in\{0,1,2\}$ where
\begin{align*}
& \hat \gamma_k=
  \begin{cases}
    0~\textrm{or}~1 &\text{if $\beta_k=1$}\\
    2&\text{if $\beta_k=0$}
  \end{cases}
\end{align*}
with given dropout probability $\eta \in [0,1]$ for the binary process $\{\beta_k: k \geq 0\}$, i.e., $\mathbb{P}(\beta_k=0)= \eta $ for all $k \geq 0$. In the case $\beta_k=1$, a transmission error may occur, independent of all other random processes, with probability $\delta \in [0,1]$. We may model the erroneous feedback channel as a discrete memoryless erasure channel with errors  depicted by a transition probability matrix
\begin{align}
\label{A_matrix}
& \mathbb{A} = (a_{ij}) = \left[\begin{array}{ccc}
    (1-\delta)(1-\eta) & \delta(1-\eta) & \eta  \\
     \delta(1-\eta) & (1-\delta)(1-\eta) & \eta \\
  \end{array} \right]
\end{align}
where $a_{ij} : = \Prob(\hat \gamma = j-1 | \gamma= i-1)$ for $i \in \{1,2\}$ and $j \in \{1,2,3\}$.
The present situation encompasses, as special cases, situations where no
acknowledgments are available (\emph{UDP-case}) and also cases where
acknowledgments are always available (\emph{TCP-case}), see also\cite{ImerYukselBasar} for a discussion in the context of closed loop control with packet dropouts. The case of perfect packet receipt acknowledgments is a special case when $\eta$ and $\delta$ above are set to zero. 

\section{Analysis of the System Model} \label{sec:aug}

\subsection{Augmented State Space Model at the Receiver} \label{sec:augsub}

To analyze the model considered in this paper, we write the dynamics of the augmented state $\theta_{k}: =[x_k ~ \hat x_{k|k-1}^s]^T$ which we want to estimate  at the remote receiver as
\begin{align*}
  & \theta_{k+1}= \mathcal{A} \theta_k + \xi_k
\end{align*}
where $\mathcal{A}:= \left[
                         \begin{array}{cc}
                           A & 0 \\
                         K_sC & A-K_sC
                         \end{array}
                       \right]$, and $\xi_k :=   \left[     \begin{array}{c}
                           w_k \\
                         K_s v_k
                         \end{array}
                       \right]$ by (\ref{ProcDyn}), (\ref{ProcObs}) and (\ref{KFR:Pred}). From (\ref{sensor_packet}), the observation  is given by $z_k = \nu_k (\hat x_{k|k}^s + q_k^x) + (1-\nu_k) (\epsilon_k^s +q_k^\epsilon)$, or
\begin{align*}  
  & z_k = \mathcal{C}(\nu_k)\theta_k + \zeta_k
\end{align*}
where $C(\nu_k) := [K_fC ~ \nu_k I-K_fC]$, and $\zeta_k := K_f v_k + v_k q_k^x +(1-v_k)q_k^\epsilon$ by (\ref{ProcObs}), (\ref{KFR:Filter}) and (\ref{InnDef}) (note that $K_fC$ is a square matrix). 
We note that $\{\xi_k: k \geq 0\}$ and $\{\zeta_k: k \geq 0\}$ are zero-mean noise processes. The covariance of the process $\{\xi_k: k \geq 0\}$ is
\begin{align*}
& Q := E[\xi_k\xi_k^T] = \left[
                         \begin{array}{cc}
                           \Sigma_w & 0 \\
                         0 & K_s \Sigma_v K_s^T
                         \end{array}
                       \right] \geq 0
\end{align*}
while the covariance of the process $\{\zeta_k: k \geq 0\}$ is given by
\begin{align*}
 & R(\nu_k) \! := \! E[\zeta_k \zeta_k^T] = K_f \Sigma_v K_f^T + \nu_k^2 \Sigma_q^x + (1-\nu_k)^2 \Sigma_q^\epsilon \geq 0.
\end{align*}
The matrix $S$ which models the correlation between the augmented state process noise $\{\xi_k: k \geq 0\}$ and the measurement noise $\{\zeta_k: k \geq 0\}$ is given by
\begin{align*}
& S := E[\xi_k\zeta_k^T] = \left[
                         \begin{array}{c}
                          0 \\
                         K_s \Sigma_v K_f^T
                         \end{array}
                       \right] .
\end{align*}
\subsection{Kalman Filter at the Receiver} \label{sec:receiver}
We assume that the receiver knows whether dropouts occurred or not, and at instances where sensor packets are received  the decision variable $\nu_k$ is also known. Therefore, the information at the receiver at time $k$, $\mathcal{Y}_k^r$, is given by the $\sigma$-field $\sigma\{\gamma_t, \gamma_t\nu_t,\gamma_t z_t: 0 \leq t \leq k\}$. We use the convention $\mathcal Y_0^r :=\{\emptyset,\Omega\}$.
At any instant $k$, the receiver estimates the process state $x_k$ through estimation of the augmented state $\theta_k$ based on the information  $ \mathcal{Y}_{k-1}^r$. We denote the conditional expectation and the associated estimation error covariance of the augmented state
\footnote{Note that if the quantization noise distribution departs from the assumed Gaussianity (Remark 
\ref{gaussquant}), then the filter at the receiver  should be interpreted as the best linear filter and 
$ \hat{\theta}_k, {\bf P}_k$ will represent the corresponding estimate and its covariance, and will only be an approximation for the conditional mean and error covariance.} as $ \hat{\theta}_k := \E[\theta_k \,|\, \mathcal{Y}_{k-1}^r]$ and
\begin{align}
\label{partitioned_matrix}
& \!\! \!\! \! {\bf P}_k \!\! := \! \E[(\theta_k-\hat{\theta}_k)(\theta_k-\hat{\theta}_k)^T|\mathcal{Y}_{k-1}^r] \! = \!\! \left[
                         \begin{array}{cc}
                           P_k^{1,1} & P_k^{1,2}  \vspace{0.1cm} \\
                         P_k^{1,2} & P_k^{2,2}
                         \end{array}
                       \right].
\end{align}
Let $\hat x^r_k := \E[x_k| \mathcal{Y}_{k-1}^r]$. Then
\begin{align*}
&  P_k^{1,1} \equiv \E[(x_k-\hat x^r_k)(x_k-\hat x^r_k)^T|\mathcal{Y}_{k-1}^r]
\end{align*}
is the state estimation error covariance at the receiver at time $k$.
The estimation error covariance ${\bf P}_{(\cdot)}$ satisfies the following random Riccati equation of Kalman filtering with correlated process and measurement noises:
\begin{align} 
& {\bf P}_{k+1} = \mathcal{A} {\bf P}_k\mathcal{A}^T + Q - \gamma_k [\mathcal{A} {\bf P}_k \mathcal{C}^T(\nu_k)+S]  \notag \\
 & \hspace{0.6cm} [\mathcal{C}(\nu_k) {\bf P}_k \mathcal{C}^T(\nu_k) + R(\nu_k)]^{-1}[\mathcal{A} {\bf P}_k \mathcal{C}^T(\nu_k)+S]^T. \label{Receiver:RiccatiEq}
\end{align}
Note that $\gamma_k$ appears as a random coefficient in the Riccati equation (\ref{Receiver:RiccatiEq}). 

\begin{thm} \label{EEC:thm}The estimation error covariance ${\bf P}_{(\cdot)}$ of the augmented system is of the form 
\begin{align}
\label{covariance_matrix_form}
& {\bf P}_k= \left[
                         \begin{array}{cc}
                           P_k^{1,1} & P_k^{1,1}-P_s \vspace{0.1cm}\\
                          P_k^{1,1}-P_s & P_k^{1,1}-P_s
                         \end{array}
                       \right], \quad k \geq 0.
\end{align}
\end{thm}
\textit{Proof}: See Appendix A.

Theorem \ref{EEC:thm} is useful in numerical solutions of the stochastic control problems considered in the next section, in that it reduces the size of the state space in which we need to consider.

\section{The Optimal Transmission Policy Problem}\label{sec:opt}
Based on the discussion in Section \ref{sec:wireless-sensor-Kalman}, the decision of whether to send the innovation $\epsilon_k^s$, i.e. set $\nu_k=0$, or the state estimate $\hat x_{k|k}^s$, i.e. set $\nu_k=1$, will result in bit rates $n_0 \equiv n(\nu_k=0)$ or $n_1 \equiv n(\nu_k=1)$, respectively, such that $n_0 \leq n_1$. To maintain a fixed packet loss probability $p$, these bit rates yield different bit error probabilities $p_{b}^0$ and $p_{b}^1$ where
\begin{align*}
& p_{b}^0 = 1-(1-p)^{1/n_0} \geq p_{b}^1 = 1-(1-p)^{1/n_1} 
\end{align*}
by (\ref{PLP}) and the fact that $n_0 \leq n_1$. The required transmission energy for bit error probabilities $p_{b}^0$ and $p_{b}^1$ will be denoted by $E_{b}^0$ and $E_{b}^1$, respectively. Since the transmission energy is a decreasing function of the bit error probability we have $E_{b}^0 \leq E_{b}^1$. For example, in the case of AWGN channel with BPSK modulation, (\ref{BEP}) implies that
\begin{align*}
& E_{b}^0 = N_0 \times \big(\textrm{erfc}^{-1}(2 p_{b}^0)\big)^2,~ E_{b}^1 = N_0 \times \big(\textrm{erfc}^{-1}(2 p_{b}^1)\big)^2
\end{align*}
where $\textrm{erfc}^{-1}(.)$ is the inverse complementary error function, which is monotonically decreasing. 

We  define the energy per packet of $n$ bits at time $k$ as $J (\nu_k) = n_{\nu_k} \times E_{b}^{\nu_k}$ which depends on the control variable $\nu_k \in\{0,1\}$.

We now aim to design optimal transmission policies in order to minimize a convex combination of the trace of the receiver's expected estimation error variance and the amount of energy required at the sensor for sending the packet to the receiver. This optimization problem is formulated as a long term average (infinite-time horizon) stochastic control problem
\begin{align}
 \label{opt_prob}
\! \min_{\{\nu_k\}} \lim\!\sup_{\!\!\!\!\!T \rightarrow \infty} \!\frac{1}{T} \!\!\sum_{k=0}^{T-1} \!\mathbb{E} \big[\lambda \textrm{Tr} P_{k\!+\!1}^{1,1} \!+\!(1\!-\!\lambda) J(\nu_k) \big|\{\hat \gamma_l\}_{0}^{\!k\!-\!1}\!,\!\{\nu_l\}_{0}^k, \!P_{x_0}\!\big] 
\end{align}
where $\lambda \in [0,1]$ is the weight, and $P_{k+1}^{1,1} $ is the submatrix of ${\bf P}_{k+1}$ (see (\ref{partitioned_matrix})) obtained from the Riccati equation (\ref{Receiver:RiccatiEq}). To take into account the fact that acknowledgements are unreliable, the expectation in (\ref{opt_prob}) is conditioned on the transmission success process of the feedback channel $\{\hat \gamma_l\}$ instead of the packet loss acknowledgment process of the forward channel $\{\gamma_l\}_{l=0}$. Thus, in problem (\ref{opt_prob}), $\nu_k$ can only depend on $\{\hat \gamma_l\}_{0}^{\!k\!-\!1}\!$, $\!\{\nu_l\}_{0}^k$, and $\!P_{x_0}$. Therefore, this formulation falls within the general framework of stochastic control problems with imperfect state information. 


\subsection{The Case of Perfect Packet Receipt Acknowledgments} \label{sec:per} First, let us assume that the smart sensor has perfect information about whether the packets have been received at the remote estimator or not, i.e. $\eta$ and $\delta$ are set to zero in Section \ref{ternary_subsec}. The optimization problem (\ref{opt_prob}) is then reduced to a stochastic control problem with perfect state information
\begin{align*}
\! \min_{\{\nu_k\}} \lim\!\sup_{\!\!\!\!\!T \rightarrow \infty}\!\frac{1}{T} \!\!\sum_{k=0}^{T-1} \!\mathbb{E} \big[\lambda \textrm{Tr} P_{k\!+\!1}^{1,1} \!+\!(1\!-\!\lambda) J(\nu_k)\big|\{\gamma_l\}_{0}^{\!k\!-\!1}\!,\!\{\nu_l\}_{0}^k, \!P_{x_0}\!\big] 
\end{align*}
which may be written as
\begin{align}
& \!\!\!\!\!\min_{\{\nu_k\}}  \lim\!\sup_{\!\!\!\! T \rightarrow \infty} \! \frac{1}{T} \sum_{k=0}^{T-1} \mathbb{E} \big[\lambda \textrm{Tr} P_{k+1}^{1,1} +(1-\lambda) J(\nu_k) \big| {\bf P}_k, \nu_k \big] \label{opt_prob:pf}
\end{align}
due to the fact that ${\bf P}_k$ is a deterministic function of $\{\gamma_l\}_{l=0}^{k-1}$, $\{\nu_l\}_{l=0}^{k-1}$, and $P_{x_0}$.
Denote
\begin{align}
& \! \! \!  \mathcal L ({\bf P},\gamma,\nu) \!  := \mathcal{A} {\bf P} \mathcal{A}^T + Q - \gamma [\mathcal{A} {\bf P} \mathcal{C}^T(\nu)+S]  \notag \\
 & \hspace{2cm} [\mathcal{C}(\nu) {\bf P} \mathcal{C}^T(\nu) + R(\nu)]^{-1}[\mathcal{A} {\bf P} \mathcal{C}^T(\nu)+S]^T \notag \\
& \hspace{1cm} \equiv  \left[
                         \begin{array}{cc}
                           \mathcal L^{1,1} ({\bf P},\gamma,\nu)&  \mathcal L^{1,1} ({\bf P},\gamma,\nu)-P_s \\
                            \mathcal L^{1,1} ({\bf P},\gamma,\nu)-P_s &   \mathcal L^{1,1} ({\bf P},\gamma,\nu)-P_s
                         \end{array}
                       \right] 
\label{RiccatiEqOperator}
\end{align}
as the random Riccati equation operator (see Theorem \ref{EEC:thm}), where matrices $\mathcal{A}$, $Q$, $\mathcal{C}$, $S$ and $R$ are given in Section \ref{sec:augsub}. 

\begin{thm}[Perfect Packet Receipt Acknowledgments]  Independent of the initial estimation error variance $P_{x_0}$, the value of problem (\ref{opt_prob:pf}) is given by $\rho$, which is the solution of the average cost optimality (Bellman) equation
\begin{align}
& \rho + V({\bf P}) = \!\!\! \min_{\nu \in \{0,1\}} \!\! \Big(\mathbb{E} \big[\lambda\textrm{Tr} \mathcal L^{1,1} ({\bf P},\gamma,\nu) + (1-\lambda) J(\nu) \big|{\bf P}, \nu\big] \notag \\
& \hspace{3.5cm} +  \mathbb{E} \big[V \big(\mathcal L ({\bf P},\gamma,\nu)\big) | {\bf P}, \nu \big] \Big) \label{perfecfeed_dp}
\end{align}
	where										 $V$ is called the relative value function.
\end{thm}

\textit{Proof}: The proof follows from the dynamic programming principle for average cost stochastic control problems (see e.g. Proposition 7.4.1 in \cite{bertsekas1995dynamic}). \qed

The stationary solution to the problem (\ref{opt_prob:pf}) is then given by
\begin{align}
& \nu^o({\bf P}) = \arg \min_{\nu \in \{0,1\}} \Big(\mathbb{E} \big[\lambda\textrm{Tr} \mathcal L^{1,1} ({\bf P},\gamma,\nu) \notag \\
& \hspace{1cm} + (1-\lambda) J(\nu) \big|{\bf P}, \nu\big]   \mathbb{E} \big[V \big(\mathcal L ({\bf P},\gamma,\nu)\big) | {\bf P}, \nu \big] \Big) \label{Optimal:Sol}
\end{align}
where $V(\cdot)$ is the solution to (\ref{perfecfeed_dp}).

\begin{remark} Equation (\ref{perfecfeed_dp}) together with the control policy $\nu^o$ defined in (\ref{Optimal:Sol}) is known as the average cost optimality equations. If a control $\nu^o$, a measurable function $V$, and a constant $\rho$ exist which solve equations (\ref{perfecfeed_dp})-(\ref{Optimal:Sol}), then the strategy $\nu^o$ is optimal, and $\rho$ is the optimal cost in the sense that
\begin{equation*}
\begin{split}
& \!\! \lim\!\sup_{\!\!\!\! T \rightarrow \infty} \! \frac{1}{T} \sum_{k=0}^{T-1} \mathbb{E} \big[\lambda \textrm{Tr} P_{k+1}^{1,1} +(1-\lambda) J(\nu_k) \big| \nu_k = \nu^o({\bf P}_k) \big] \!\!= \! \rho \\
\end{split}
\end{equation*}
and for any other control policy $\{\nu_k \in \{0.1\}: k \geq 0\}$, \begin{equation*}
\begin{split}
& \!\! \lim\!\sup_{\!\!\!\! T \rightarrow \infty} \! \frac{1}{T} \sum_{k=0}^{T-1} \mathbb{E} \big[\lambda \textrm{Tr} P_{k+1}^{1,1} +(1-\lambda) J(\nu_k) \big| \nu_k \big] \! \geq \! \rho \\
\end{split}
\end{equation*}
The reader is referred to\cite{arapostathis1993discrete} for a proof of the average cost optimality equations and related results. We solve the Bellman equation by the use of {\it relative value iteration} algorithm (see Chapter 7 in \cite{bertsekas1995dynamic}). 
\end{remark} 

In (\ref{perfecfeed_dp}), the term $ \E \big[\mathcal L^{1,1} ({\bf P},\gamma,\nu) | {\bf P}, \nu] $ is the submatrix (similar to (\ref{partitioned_matrix})) of the following matrix
\begin{align}
& \!\!\!  \E \big[\mathcal L ({\bf P},\gamma,\nu) | {\bf P}, \nu]\! = \! A {\bf P} A^T \!\! + Q \! - \! (1-p) \! \times \! [A {\bf P} C^T(\nu)+S] \notag \\ & \hspace{1cm} \times [C(\nu) {\bf P} C^T(\nu) + R(\nu)]^{-1}[A {\bf P} C^T(\nu)+S]^T  \label{ICFE}
\end{align}
where $p$ is the packet loss probability of the forward erasure communication channel given in Section \ref{sec:comm-chann}.

\subsection{The Case of Imperfect Packet Receipt Acknowledgments}\label{sec:Imperfeed}

In the formulation of problem (\ref{opt_prob}), the smart sensor does not have perfect knowledge about whether its transmissions have been received at the receiver. Hence, at time $k$ the sensor has only ``imperfect state information'' about $\{{\bf P}_t: 1 \leq t \leq k\}$ via the acknowledgment process $\{\hat \gamma_t, 0 \leq t \leq k-1 \}$.
We will reduce the optimization problem (\ref{opt_prob}) to a stochastic control problem with perfect state information by using the notion of information-state \cite{kumar1986stochastic}. 
For $k \geq 0$ denote $z^k := \{\hat \gamma_0, \cdots, \hat \gamma_{k}, \nu_0, \cdots, \nu_{k-1},P_{x_0}\}$ as all observations about the receiver's Kalman filtering state estimation error covariance at the sensor after  transmission at time $k$ and before transmission at time $k+1$. We set $z^{-1}:=\{P_{x_0}\}$. The {\it information-state} is defined by
\begin{align}
& f_{k+1}({\bf P}_{k+1}|z^k,\nu_k) = \mathbb{P}({\bf P}_{k+1}|z^k,\nu_k), \qquad k \geq 0
\end{align}
which is the conditional probability of the estimation error covariance ${\bf P}_{k+1}$ given $(z^k,\nu_k)$. The following lemma shows how $f_{k+1}(\cdot|z^k,\nu_k)$ can be determined from $f_k(\cdot|z^{k-1},\nu_{k-1})$ together with $\hat \gamma_k$ and $\nu_k$.

\begin{lem} The information-state $f_{(\cdot)}$ satisfies the recursion:
\begin{align} 
& f_{k+1}({\bf P}_{k+1}|z^k,\nu_k) =  \sum_{\gamma_k \in\{0,1\}}  \!\!\! \Big[ \int_{{\bf P}_k} \Big(\mathbb{P}({\bf P}_{k+1} | {\bf P}_k,\gamma_k,\nu_k) \notag \\
&   \times f_{k} ({\bf P}_k|z^{k-1},\nu_{k-1})\Big)d{\bf P}_k \times   \frac{\mathbb{P}(\hat \gamma_k|\gamma_k) \times \mathbb{P}(\gamma_k)}{ \sum_{\gamma_k \in \{0,1\}}\mathbb{P}(\hat \gamma_k|\gamma_k) \times \mathbb{P}(\gamma_k)} \Big] \notag \\
 & =: \Phi \big[f_k(\cdot|z^{k-1},\nu_{k-1}),\hat \gamma_k, \nu_k \big]({\bf P}_{k+1}).  \quad k \geq 0 \label{InfState_recur_tot}
\end{align}
with $f_0({\bf P}_0|z^{-1}) = \delta({\bf P}_0)$, where $\delta$ is the Dirac delta function. 
\label{lemmainfstaterec}
\end{lem}
\textit{Proof}: See Appendix B.

Note that in (\ref{InfState_recur_tot}) the probabilities $\mathbb{P}(\hat{\gamma}_k|\gamma_k)$ can be obtained from the probability transition matrix $\mathbb{A}$ in (\ref{A_matrix}). 
It is important to note that $\Phi$ in (\ref{InfState_recur_tot}) depends on the entire function $f_{k}(\cdot|z^{k-1},\nu_{k-1})$ and not just its value at any particular ${\bf P}_{k}$. 

We now reduce  problem  (\ref{opt_prob}) to a problem with perfect state information, where its state is given by the information state $f_{(\cdot)}$ which evolves based on the recursion (\ref{InfState_recur_tot}). 
Define the class of  matrices $\mathbb{S}$ as
\begin{align} \label{Class:S}
& \mathbb{S} := \left\{{\bf P}=\left[
                         \begin{array}{cc}
                           P &  P-P_s \\
                          P-P_s & P-P_s
                         \end{array}
                       \right]: P \geq P_s \right\}.
\end{align} 

\begin{thm}[Imperfect Packet Receipt Acknowledgments]  Independent of the initial estimation error variance $P_{x_0}$, the value of  problem (\ref{opt_prob}) is given by $\rho$, which is the solution of the average cost optimality (Bellman) equation
\begin{align}
& \rho + V(\pi) = \!\!\! \min_{\nu \in \{0,1\}} \! \Big(\mathbb{E} \big[\lambda\textrm{Tr}\mathcal L^{1,1} ({\bf P},\gamma,\nu) + (1-\lambda) J(\nu) \big|\pi, \nu\big] \notag \\
& \hspace{3.5cm} +  \mathbb{E} \big[V \big(\Phi\big(\pi,\hat \gamma, \nu)\big) | \pi, \nu \big] \Big) \label{Imperfecfeed_dp}
\end{align}
for $\pi \in \Pi$, where the operator $\Phi$ is defined in (\ref{InfState_recur_tot}), $V$ is  the relative value function, and $\Pi$ is the space of probability density functions on matrices $\mathbb{S}$ of the form (\ref{Class:S}).
\end{thm}

{\it Proof}: The proof follows from the dynamic programming principle for stochastic control problems with imperfect state information (see Theorem 7.1 in \cite{kumar1986stochastic}). \qed

Note that in  (\ref{Imperfecfeed_dp}) the state is the entire probability density function $\pi$ which takes its values in the space of probability densities $\Pi$. 
We may write the terms in (\ref{Imperfecfeed_dp}) as
\begin{align*}
&  \mathbb{E} \big[\mathcal L({\bf P},\gamma,\nu) \big| \pi, \nu] = \int_{{\bf P}} \big( \mathcal{A} {\bf P}  \mathcal{A}^T + Q) \pi ({\bf P}) d{\bf P} \\
& \hspace{0.1cm}  - (1-p) \times \int_{{\bf P}} \Big([ \mathcal{A} {\bf P}  \mathcal{C}^T(\nu)+S] [ \mathcal{C}(\nu) {\bf P}  \mathcal{C}^T(\nu) + R(\nu)]^{-1}  \notag \\ 
& \hspace{3cm} \times [ \mathcal{A} {\bf P}  \mathcal{C}^T(\nu)+S]^T \Big) \pi ({\bf P}) d{\bf P}\notag
\end{align*}
and $\mathbb{E} \big[V \big(\Phi\big(\pi,\hat \gamma, \nu)\big) | \pi, \nu \big] = \Prob(\hat \gamma =0) V \big(\Phi\big(\pi,0, \nu) \big)$ $+ \Prob(\hat \gamma =1) V \big(\Phi\big(\pi,1, \nu) \big) + \Prob(\hat \gamma =2) V \big(\Phi\big(\pi,2, \nu) \big).$

\section{A Suboptimal Transmission Policy Problem}
\label{sec:subopt}

To obtain the optimal transmission strategy in the case of imperfect packet receipt acknowledgments presented in Section \ref{sec:Imperfeed} we need to compute the solution of the Bellman equation (\ref{Imperfecfeed_dp}) in the space of probability density functions $\Pi$, which is computationally demanding. In this section we consider suboptimal policies which are computationally much less intensive than finding the optimal solution. 

We formulate the suboptimal optimization problem as
\begin{align}
&  \!\! \min_{\{\nu_k\}} \lim\!\sup_{\!\!\!\! T \rightarrow \infty} \! \frac{1}{T} \sum_{k=0}^{T-1} \mathbb{E} \big[\lambda \textrm{Tr} \hat P_{k+1}^{1,1} +(1-\lambda) J(\nu_k)  \big| {\bf \hat P}_k, \nu_k \big] \label{Imperfecfeed_prob:sub}
\end{align}
where $\hat P_{(\cdot)}^{1,1}$ is the submatrix (similar to (\ref{partitioned_matrix})) of $\hat {\bf P}_{(\cdot)}$ which is an estimate of ${\bf P}_{(\cdot)}$ computed by the  sensor based on the following recursive equations (with $\hat {\bf P}_0  = {\bf P}_0$): 

(i) In the case $\hat \gamma_k = 0$ we have
\begin{align*}
&\hat {\bf P}_{k+1} \!:=\!\! \Big(\mathcal{A} \hat {\bf P}_k  \mathcal{A}^T + Q \Big) \!\! \times\!\! \frac{\Prob(\hat \gamma_k=0|\gamma_k=0) \times \Prob(\gamma_k=0)}{ \sum_{\gamma_k \in \{0,1\}}\Prob(\hat \gamma_k=0|\gamma_k) \times \Prob(\gamma_k)} \\
& \hspace{0.2cm} + \Big( \mathcal{A} \hat {\bf P}_k  \mathcal{A}^T + Q - [ \mathcal{A} \hat {\bf P}_k  \mathcal{C}^T(\nu_k)+S]  \notag \\
& \hspace{0.7cm}  \times [ \mathcal{C}(\nu_k) \hat {\bf P}_k C^T(\nu_k) + R(\nu_k)]^{-1}[ \mathcal{A} \hat {\bf P}_k  \mathcal{C}^T(\nu_k)+S]^T \Big) \\
& \hspace{0.2cm} \times \frac{\Prob(\hat \gamma_k=0|\gamma_k=1) \times \Prob(\gamma_k=1)}{ \sum_{\gamma_k \in \{0,1\}}\Prob(\hat \gamma_k=0|\gamma_k) \times \Prob(\gamma_k)}. 
\end{align*}

(ii) in the case $\hat \gamma_k = 1$ we have
\begin{align*}
& \hat P_{k+1} \!:=\!\! \Big(\mathcal{A} \hat {\bf P}_k  \mathcal{A}^T + Q \Big)  \!\! \times\!\!  \frac{\Prob(\hat \gamma_k=1|\gamma_k=0) \times \Prob(\gamma_k=0)}{ \sum_{\gamma_k \in \{0,1\}}\Prob(\hat \gamma_k=1|\gamma_k) \times \Prob(\gamma_k)} \\
& \hspace{0.2cm} + \Big( \mathcal{A} \hat {\bf P}_k  \mathcal{A}^T + Q - [ \mathcal{A} \hat {\bf P}_k  \mathcal{C}^T(\nu_k)+S]  \notag \\
& \hspace{0.7cm}  \times [ \mathcal{C}(\nu_k) \hat {\bf P}_k C^T(\nu_k) + R(\nu_k)]^{-1}[ \mathcal{A} \hat {\bf P}_k  \mathcal{C}^T(\nu_k)+S]^T \Big) \\
& \hspace{0.2cm} \times \frac{\Prob(\hat \gamma_k=1|\gamma_k=1) \times \Prob(\gamma_k=1)}{ \sum_{\gamma_k \in \{0,1\}}\Prob(\hat \gamma_k=1|\gamma_k) \times \Prob(\gamma_k)}. 
\end{align*}

(iii) In the case $\hat \gamma_k = 2$ we have
\begin{align*}
& \hat P_{k+1} := \mathcal{A} \hat {\bf P}_k  \mathcal{A}^T + Q - \Prob(\gamma_k =1) \times [ \mathcal{A} \hat {\bf P}_k \mathcal{C}^T(\nu_k)+S]  \notag \\
& \hspace{1.2cm}  \times [ \mathcal{C}(\nu_k) \hat {\bf P}_k  \mathcal{C}^T(\nu_k) + R(\nu_k)]^{-1}[ \mathcal{A} \hat {\bf P}_k  \mathcal{C}^T(\nu_k)+S]^T. 
\end{align*}

The reason that the solution to the stochastic control problem (\ref{Imperfecfeed_prob:sub}) is only suboptimal is that the true error covariance matrix ${\bf P}_{(\cdot)}$ in (\ref{opt_prob}) is replaced by its estimate $\hat {\bf P}_{(\cdot)}$ in (\ref{Imperfecfeed_prob:sub}). The intuition behind these recursive equations can be explained as follows. Note that in the case of perfect feedback acknowledgements (\ref{Receiver:RiccatiEq}), the error covariance is updated as ${\bf P}_{k+1} = \mathcal{A} {\bf P}_k\mathcal{A}^T + Q $ in case $\gamma_k=0$, and ${\bf P}_{k+1} = \mathcal{A} {\bf P}_k\mathcal{A}^T + Q - \gamma_k [\mathcal{A} {\bf P}_k \mathcal{C}^T(\nu_k)+S] \times [\mathcal{C}(\nu_k) {\bf P}_k \mathcal{C}^T(\nu_k) + R(\nu_k)]^{-1}\times [\mathcal{A} {\bf P}_k \mathcal{C}^T(\nu_k)+S]^T$ in case $\gamma_k=1$. In our imperfect acknowledgement model, even when it is received, errors can occur such that $\hat{\gamma}_k=0$ is received when $\gamma_k=1$, and $\hat{\gamma}_k=1$ is received when $\gamma_k=0$. Thus the recursions given in (i) and (ii) are the weighted combinations of the error covariance recursions (based on the Bayes' rule using corresponding error event probabilities) in the case of perfect feedback acknowledgements. In the case $\hat{\gamma}_k=2$ where an erasure 
occurs, taking the average of the error covariances in the cases $\gamma_k=0$ and $\gamma_k=1$ is intuitively a reasonable thing to do, which motivates the recursion in (iii).

Note that $\mathbb{P}(\hat \gamma_k)= \sum_{\gamma_k \in \{0,1\}} \Prob(\hat \gamma_k | \gamma_k) \Prob(\gamma_k)$, where the conditional probabilities are given in Section \ref{ternary_subsec}. This together with the recursive equations of $\hat {\bf P}_{(\cdot)}$ implies that the expression $\mathbb{E} [\hat P^{1,1}_{k+1}| \hat {\bf P}_k,  \nu_k]$ is of the same form as $\mathbb{E} [P^{1,1}_{k+1}| {\bf P}_k, \nu_k] $ when ${\bf P}_k$ is replaced by $\hat {\bf P}_k$, and the Bellman equation for problem (\ref{Imperfecfeed_prob:sub}) is given by a similar Bellman equation to (\ref{perfecfeed_dp}). The details are omitted for brevity.
\vspace*{-0.5cm}
\section{Structural Results On the Optimal Transmission Policies for Scalar Systems} \label{sec:thresh}
This section presents structural results of the optimal transmission policies for scalar systems (where we will set $A=a$, $C=1$, $\Sigma_w=\sigma_w^2$, $\Sigma_v=\sigma_v^2$, $\Sigma_q = \sigma_q^2$) in the perfect packet receipt acknowledgments case examined in Section \ref{sec:per}, which is also valid for the suboptimal solution presented in Section \ref{sec:subopt}. The idea is to apply the submodularity concept (see \cite{ngo2009optimality,topkis2001supermodularity}) to the recursive Bellman equation (\ref{perfecfeed_dp}),
to show that the optimal policy
$\nu^o(\cdot)$ in both scenarios is monotonically increasing with respect to the receiver's state estimation error variance ${\bf P}^{1,1}$. This monotonicity then implies a threshold structure since the control space has only two elements $\{0,1\}$.

\begin{defi}[\cite{ngo2009optimality} after \cite{topkis2001supermodularity}] A function $F(x, y) : X \times Y \rightarrow S$ is submodular in $(x, y)$ if $F(x_1, y_1) + F(x_2, y_2) \leq F(x_1, y_2) + F(x_2, y_1)$ for all $x_1, x_2 \in X$ and $y_1, y_2 \in Y$ such that $x_1 \geq x_2$ and $y_1 \geq y_2$. \qed
\end{defi} 

It is important to note that the submodularity is a sufficient condition for optimality of monotone increasing policies. Specifically, if $F(x, y)$ defined above is submodular in $(x, y)$ then $y(x) =\arg\min_y F(x, y)$ is non-decreasing in $x$.

We define the ordering $\geq$ for matrices in class $\mathbb{S}$ of the form (\ref{Class:S})  as
${\bf P}_1 \geq {\bf P}_2$
if ${\bf P}_1-{\bf P}_2$ is positive semi-definite. It is evident that for ${\bf P}_1, {\bf P}_2 \in \mathbb{S}$ we have ${\bf P}_1 \geq {\bf P}_2$ if and only if $P_1^{1,1} \geq P_2^{1,1}$. We also define $F: \mathbb{S} \times \{0,1\} \rightarrow \mathbb{S}$ as 
\begin{align*}
& F({\bf P},\nu) =  \mathcal{A} {\bf P}  \mathcal{A}^T + Q - (1-p) \times [ \mathcal{A} {\bf P}  \mathcal{C}^T(\nu)+S] \\ 
& \hspace{2cm} \times [ \mathcal{C}(\nu) {\bf P}  \mathcal{C}^T(\nu) + R]^{-1}[ \mathcal{A} {\bf P}  \mathcal{C}^T(\nu)+S]^T
\end{align*}
based on the instantaneous cost $ \E \big[\mathcal L ({\bf P},\gamma,\nu) | {\bf P}, \nu]$  in (\ref{ICFE}). Note that in the scalar case $R$ can be made independent of $\nu_k$. 

\begin{lem}\label{Submodul:lem} The function $F({\bf P},\nu)$ is submodular in $({\bf P},\nu)$, i.e., for ${\bf P}_1, {\bf P}_2 \in \mathbb{S}$ such that ${\bf P}_1 \geq {\bf P}_2$ we have
\begin{align}  \label{state-cost_sub}
& \!\! F^{1,1}({\bf P}_1, 1) \! + F^{1,1}({\bf P}_2, 0) \leq F^{1,1}({\bf P}_1, 0) + \!F^{1,1}({\bf P}_2, 1)
\end{align} 
where $F^{1,1}(\cdot,\cdot)$ is the (1,1) entry of $F(\cdot,\cdot)$. This implies that $F({\bf P}_1, 1) +F({\bf P}_2, 0) \leq F({\bf P}_1, 0) + F({\bf P}_2, 1).$
\end{lem}
\textit{Proof}: See Appendix C.

We now present the relative value iteration algorithm to solve the Bellman equation (\ref{perfecfeed_dp}). It is used to construct structural results for the optimal transmission policy.
First, we consider the Bellman equation for the finite $T$-horizon stochastic control problem: 
\begin{align}
& V_t({\bf P}) =  \min_{\nu \in \{0,1\}} \Big(\mathbb{E} \big[\lambda\mathcal L^{1,1} ({\bf P},\gamma,\nu) + (1-\lambda) J(\nu) \big|{\bf P}, \nu\big] \notag \\
& \hspace{0.8cm} +  \mathbb{E} \big[V_{t+1} \big(\mathcal L ({\bf P},\gamma,\nu)\big) | {\bf P}, \nu \big] \Big), \quad 0 \leq t \leq T-1\label{perfecfeed_dp:FiniteT}
\end{align}
with terminal condition $V_T({\bf P}) = 0$ where $T$ is large. We now define the function
\begin{align}
& H_t (\cdot) := V_t(\cdot) - V_t({\bf P}_f),  \quad 0 \leq t \leq T
\end{align}
where ${\bf P}_f \neq {\bf P}_0$ is fixed. We then have the following relative value iteration algorithm recursion
\begin{align}
& H_t({\bf P}) = \min_{\nu \in \{0,1\}} \Big(\mathbb{E} \big[\lambda\mathcal L^{1,1} ({\bf P},\gamma,\nu) + (1-\lambda) J(\nu) \big|{\bf P}, \nu\big] \notag \\
& \hspace{2cm} + \mathbb{E} \big[V_{t+1} \big(\mathcal L ({\bf P},\gamma,\nu)\big) | {\bf P}, \nu \big] \Big) \notag \\
& \hspace{0.5cm} - \min_{\nu \in \{0,1\}} \Big(\mathbb{E} \big[\lambda\mathcal L^{1,1} ({\bf P},\gamma,\nu) + (1-\lambda) J(\nu) \big|{\bf P}={\bf P}_0, \nu\big] \notag \\
& \hspace{2cm} + \mathbb{E} \big[V_{t+1} \big(\mathcal L ({\bf P},\gamma,\nu)\big) | {\bf P}={\bf P}_f, \nu \big] \Big)  \label{RVFI}
\end{align}
for $0 \leq t \leq T-1$. It can be shown that the relative value recursion (\ref{RVFI}) converges to the optimal solution $\rho$ of the infinite-time horizon average cost Bellman equation (\ref{perfecfeed_dp}) such that $\rho \approx H_0({\bf P}_0)$ (see the discussion on page 391 in Chapter 7 of \cite{bertsekas1995dynamic}).

\begin{thm}\label{Submodul:thm} The optimal transmission policy in the case of perfect feedback channel  is threshold with respect to the receiver's state estimation error variance $P^{1,1}$ (and hence in the augmented state estimation error covariance ${\bf P})$,
i.e., 
\begin{align} \label{THRES:Polic}
& \nu^o ({\bf P}) = \left\{ \begin{array}{cl}
0  , & \textrm{if} ~ P_k^{1,1} \leq \phi^*   \\
1  , & \textrm{otherwise}
\end{array} \right. 
\end{align}
where $\phi^*$ is the threshold.
\end{thm}

\textit{Proof}: See Appendix D.

\begin{figure}[!t]
\begin{center}
\includegraphics[width=8.5cm,height=5cm]{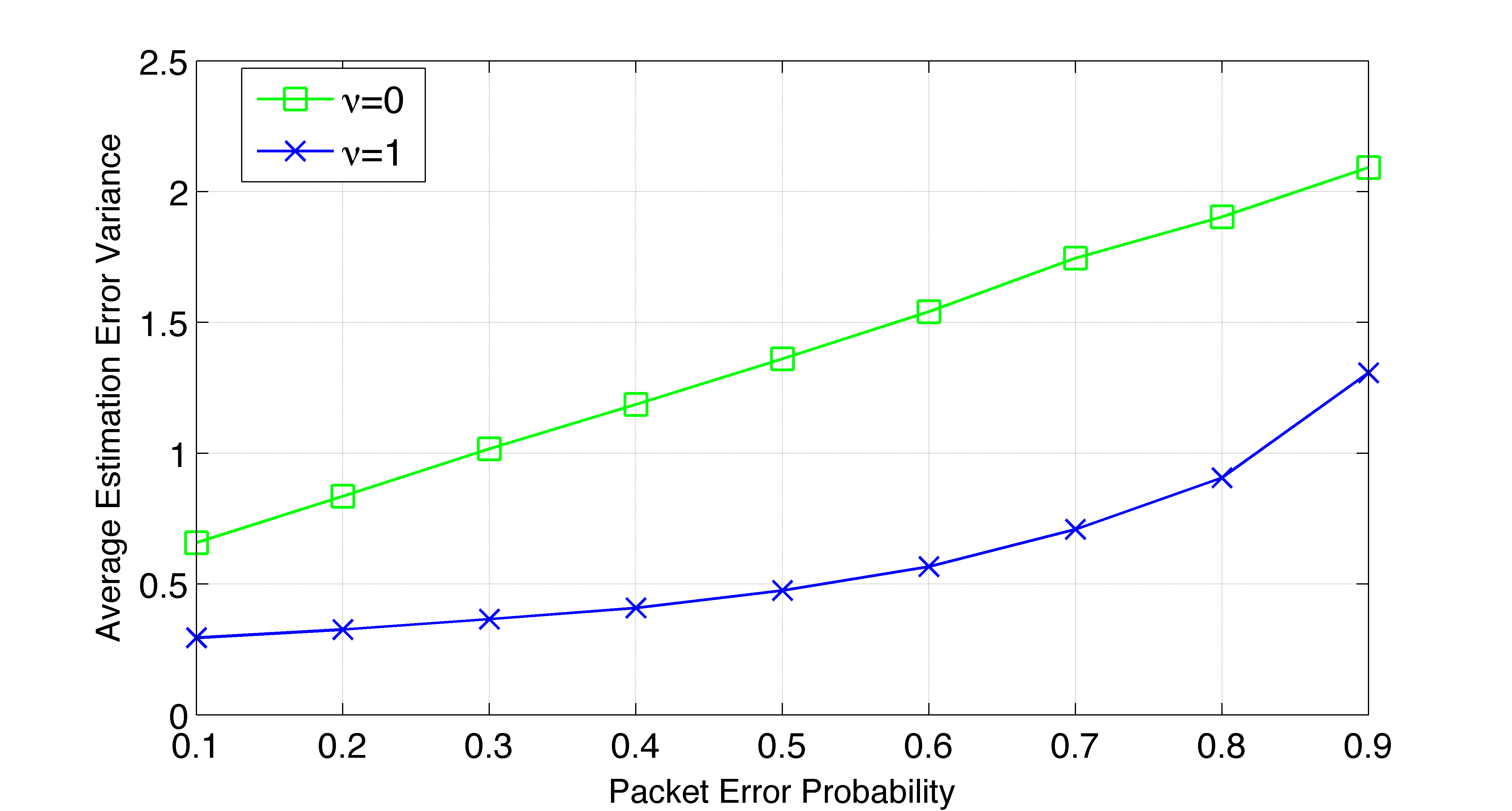}  
\caption{Perfect feedback case: Average estimation error variance versus \newline  the packet error probabilities for the two cases of $\nu=0$ and $\nu=1$} 
\label{Simu1} 
\end{center}                             
\end{figure}

The threshold structure of Theorem (\ref{Submodul:thm}) simplifies the implementation of the optimal transmission policy significantly. However, this requires knowledge of the threshold $\phi^*(\cdot)$. In general, there is no closed form expression for $\phi^*(\cdot)$, but it can be found via iterative search algorithms. Here we present a stochastic gradient  algorithm based on Algorithm 1 in \cite{krishnamurthy2012sequential}.

First, we establish some notation. For fixed ${\bf P}$ denote
\begin{align*}
& J (\theta^*) :=  \mathbb{E} \big[\lambda\mathcal L^{1,1} ({\bf P},\gamma,\nu^o) + (1-\lambda) J(\nu) \big|{\bf P}, \nu^o\big] \notag \\
& \hspace{1.5cm} +  \mathbb{E} \big[V_{0} \big(\mathcal L ({\bf P},\gamma,\nu^o)\big) | {\bf P}, \nu^o \big] 
\end{align*}
where the policy $\nu^o$ is defined in (\ref{THRES:Polic}) based on the threshold $\phi^* $, and $V_0(\cdot)$ is obtained from the finite $T$-horizon Bellman equation (\ref{perfecfeed_dp:FiniteT}). For $n \in \mathbb{N}$, $0.5 < \kappa \leq 1$ and $\omega,\varsigma >0$ we denote $\omega_n:=\frac{\omega}{(n+1)^\kappa}$ and $\varsigma_n := \frac{\varsigma}{(n+1)^\kappa}$.

{\it Stochastic gradient algorithm for computing the threshold}. For fixed ${\bf P}$ in the relative value algorithm (\ref{RVFI}) the following steps are carried out:
 
Step 1) Choose the initial threshold $\phi^{(0)}$.

Step 2) For iterations $n=0,1, \cdots$
\begin{itemize}
\item Compute the gradient:
\begin{align} \label{gra} 
& \!\! \partial_{\phi} J_n := \frac{J(\phi^{(n)}+\omega_nd_n)- J(\phi^{(n)}-\omega_nd_n)}{2 \omega_n}d_n
\end{align} 
where $d_n \in \{-1,1\}$ is a random variable such that $\Prob(d_n=-1)=\Prob(d_n=1)=0.5$.
\item Update the threshold via $\phi^{(n+1)} = \phi^{(n)} -\varsigma_n  \partial_{\phi} J_n$ which gives
\begin{align*} 
& \nu^{(n+1)} ({\bf P}) = \left\{ \begin{array}{cl}
0 , & \quad \textrm{if} ~ P^{1,1} \leq \phi^{(n+1)}  \\
1 , & \quad  \textrm{otherwise.}
\end{array} \right. 
\end{align*}
\end{itemize}
The above algorithm is a gradient-estimate 
based algorithm (see \cite{spall2005introduction}) for estimating the 
optimal threshold $\phi^*(\cdot)$ where only measurements of the loss function is available (i.e., no gradient information). We note that (\ref{gra}) evaluates an approximation to the gradient. This algorithm generates a sequence of estimates for the threshold policy $\phi^*$ which converges to a local minimum with corresponding energy allocation $\nu^*$. The reader is referred to \cite{spall2005introduction} for associated convergence analysis of this and other related algorithms (see e.g., Theorem 7.1 in \cite{spall2005introduction}). Note that gradient-estimate based algorithms are sensitive to initial conditions and should be evaluated for several distinct initial conditions to find the best local minimum.

\begin{figure}[!t]
\begin{center}
\includegraphics[width=8.5cm,height=5cm]{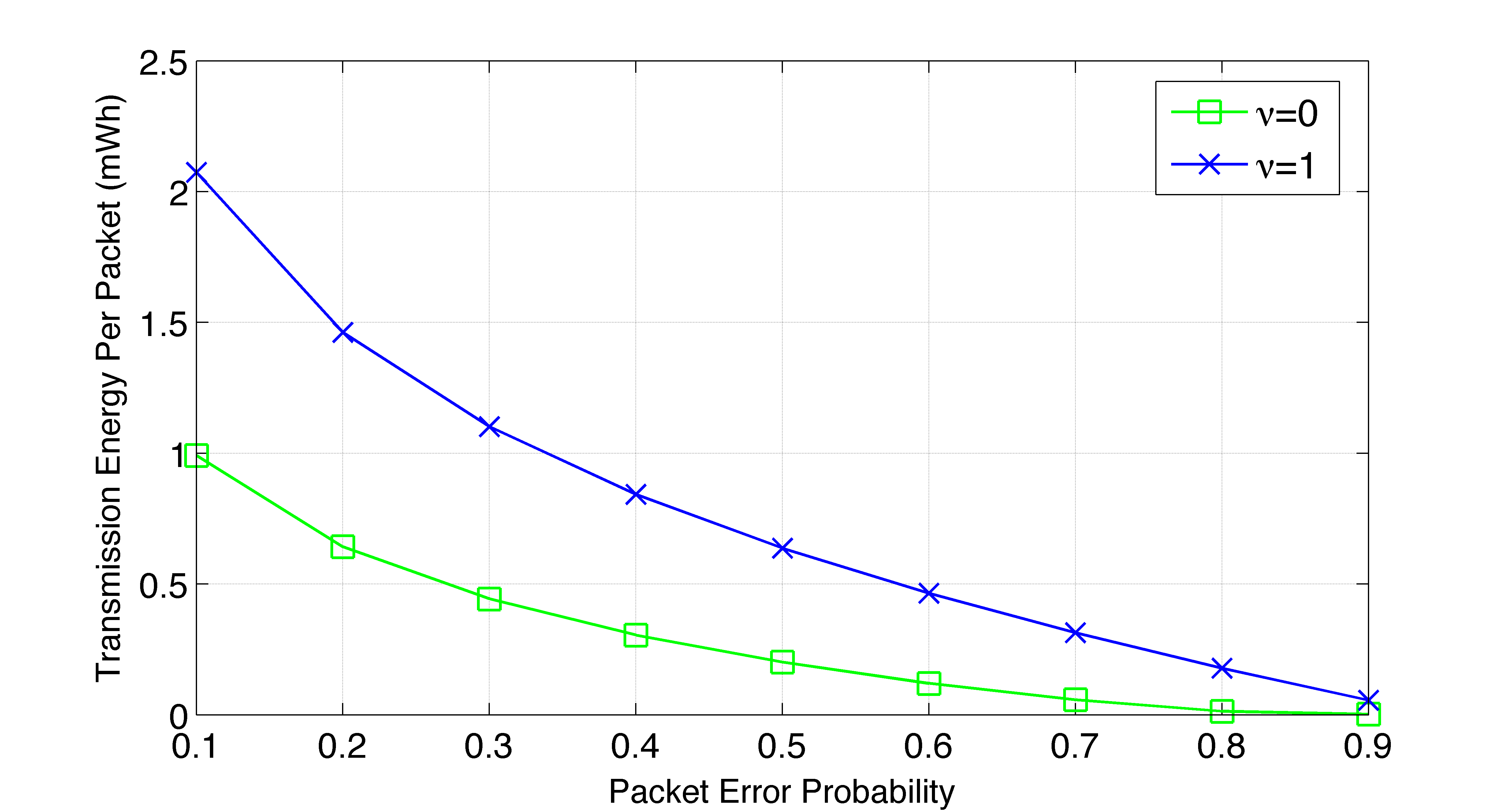}   
\caption{Perfect feedback case: Transmission energy per packet (mWh) \newline versus the packet error probabilities for the two cases of $\nu=0$ and $\nu=1$} 
\label{Simu2} 
\end{center}                                
\end{figure}

\section{Numerical Examples} \label{numerical_sec} We present here numerical results for a scalar  model with parameters $a=0.95$, $c=1$, $\sigma_w^2=0.25$, $\sigma_v^2=0.01$ and $P_{x_0}$=1 in (\ref{ProcDyn}) and (\ref{ProcObs}). These values give $P_s=0.26$, $K_s=0.91$, $K_f=0.96$ and $\Sigma_s=2.30$ in Section \ref{sec:wireless-sensor-Kalman}. We take $\sigma_q^2=0.01$ in (\ref{Quant-var}) together with an optimal Lloyd-Max quantizer which yields $n_0=3$ and $n_1=5$ by (\ref{LMQuzn}). In the simulation results, an AWGN channel with BPSK modulation is assumed where $N_0 =0.01$ in (\ref{BEP}) (see Section \ref{sec:comm-chann}).

\subsection{Perfect Feedback Communication Channel Case}
 First, let the packet error probability $p$ in (\ref{PLP}) be equal to 0.2. This gives $p_b^0=0.07$ and $p_b^1=0.04$, and hence, energy per bit levels of $E_b^0=0.21$ and $E_b^1=0.29$, see Section \ref{sec:opt}. 

In Fig \ref{Simu1}, we plot the average estimation error variance versus the packet error probabilities. More precisely, we take $\lambda=1$ in (\ref{opt_prob:pf}) without computing the optimal solution.  Instead, we let the transmission policies $\{\nu_k, k \geq 0\}$ be fixed either equal to zero (sending innovations) or one (sending state estimates). 
On the other hand, Fig \ref{Simu2} presents the packet transmission energy $J(\nu)$ (in milliwatt hour (mWh)) defined in Section \ref{sec:opt} versus the packet error probabilities. We let the transmission policy $\nu$ be fixed equal to either zero (sending innovations) or one (sending state estimates). 
Figs. \ref{Simu1} and \ref{Simu2} show that transmitting local estimates gives smaller error covariance, but also requires more transmit energy, than transmitting local innovations, which motivates the optimization formulation (\ref{opt_prob}).

\begin{figure}[!t]
\begin{center}
\includegraphics[width=8.5cm,height=5cm]{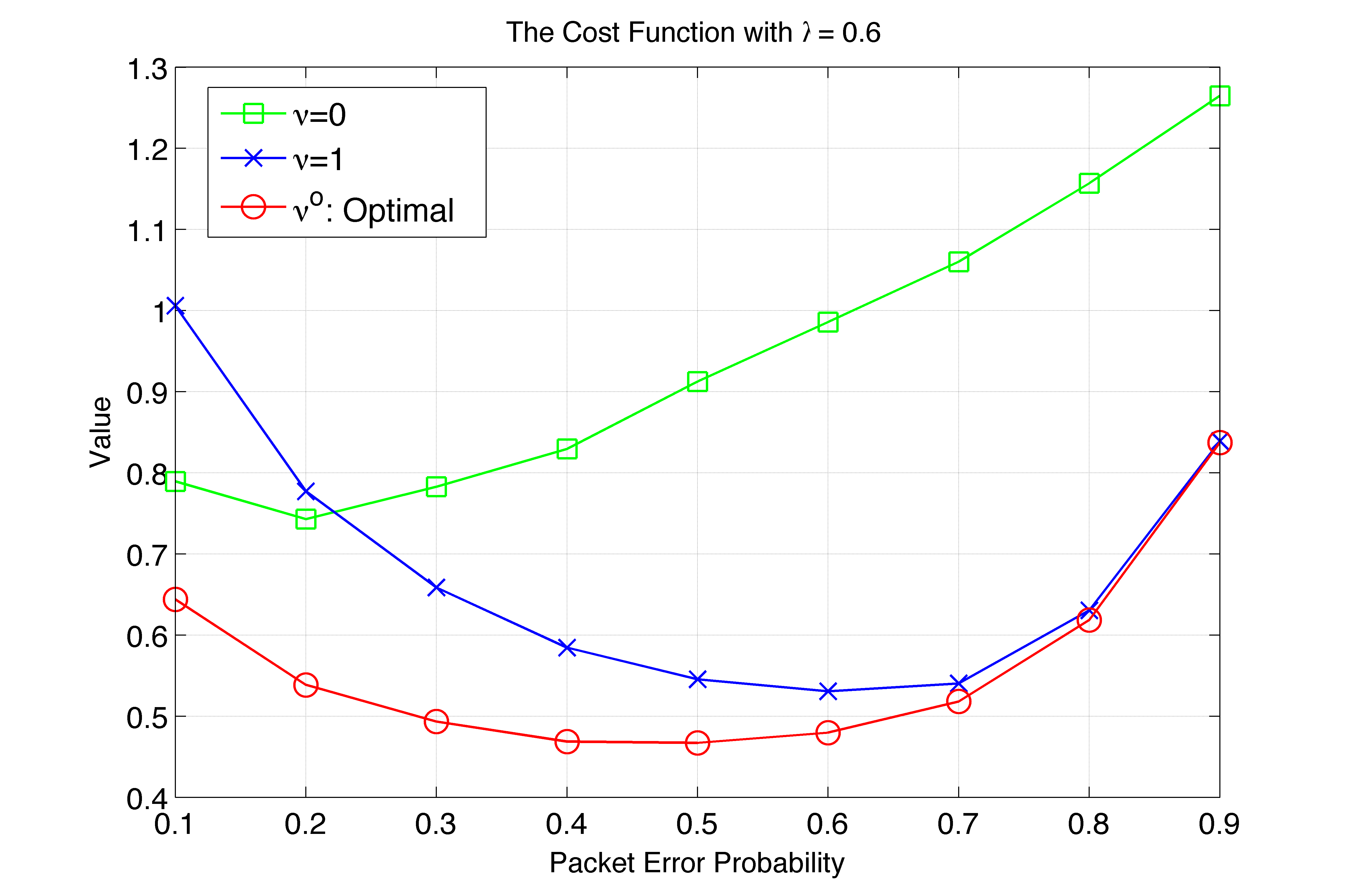}  
\caption{\footnotesize Perfect feedback case: Performance versus the packet error probabilities} \label{Simu3} 
\end{center}                                
\end{figure}

We now set the weight $\lambda$ in (\ref{opt_prob:pf}) to 0.6. The discretized equation of the relative value algorithm (\ref{RVFI}) is used for the numerical computation of the optimal transmission policy. In solving the Bellman equation (\ref{perfecfeed_dp}) we use 40 discretization points for the state estimate error variance $P_k^{1,1}$ in the range of [0,2]. In Fig. \ref{Simu3} we plot the  convex combination of the receiver's expected estimation error variance and the energy needed to transmit the packets, versus the packet loss probability $p \in [0.1,0.9]$ for the cases of: (i) fixed transmission policy $\nu=0$, (ii) fixed transmission policy $\nu=1$, and (iii) optimal transmission policy $\nu^o$. We observe that for small packet loss probabilities sending innovations ($\nu=0$) is better than sending the state estimates ($\nu=1$). On the other hand, for large packet loss probabilities sending the state estimates gives better performance than sending the innovations, due to the poor estimation performance when sending innovations when the packet loss probability is high. 


\subsubsection*{Threshold Policy} 
\begin{figure}[!t]
\begin{center}
\includegraphics[width=8.5cm,height=5cm]{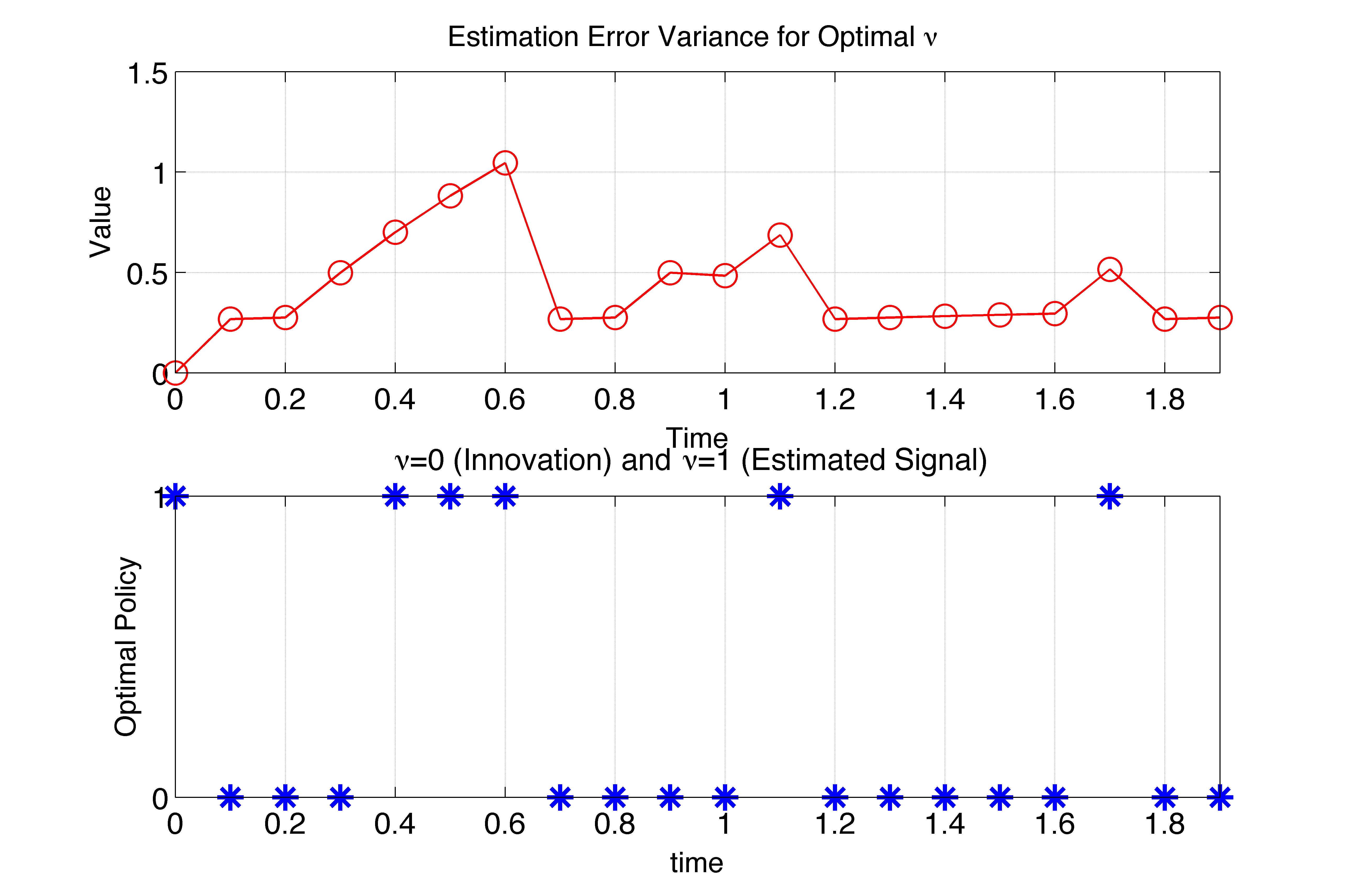}
\caption{Perfect feedback case: A single simulation run}
\label{Simu5}
\end{center}
\end{figure}

Let the packet error probability $p$ in (\ref{PLP}) be equal to 0.2. Applying the stochastic gradient algorithm given at the end of Section \ref{sec:thresh} with parameters $\omega=0.3, \varsigma =0.5$ and $\kappa=1$  yields 
the threshold $\phi^*=0.5$. 
For this case, a single run simulation result of the receiver's state estimation error variance $P_{(\cdot)}^{1,1}$ is given together with the optimal transmission strategy in Fig. \ref{Simu5}.

\begin{figure}[!t]
\begin{center}
\includegraphics[width=8.5cm,height=5cm]{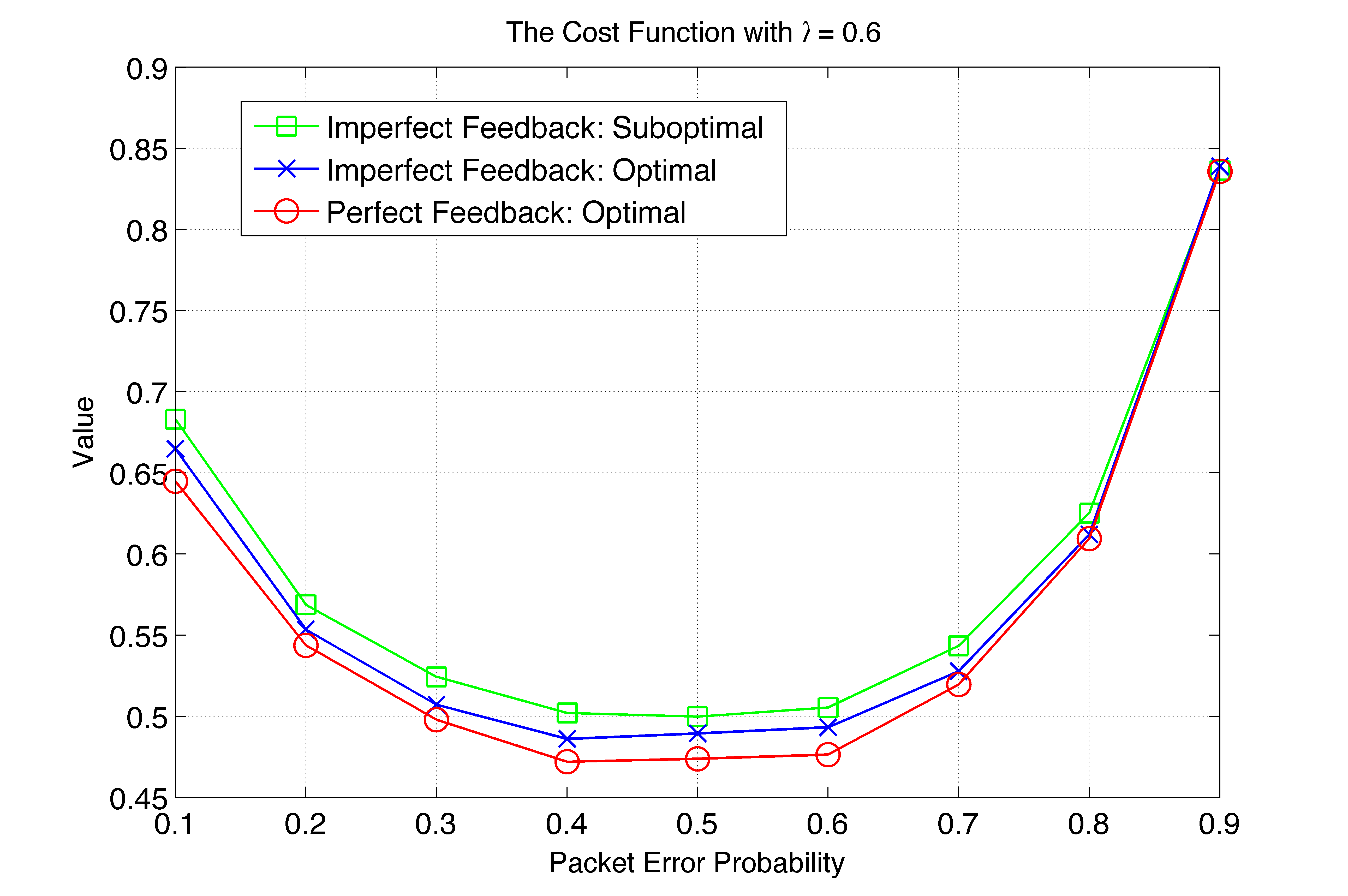}  
\caption{Performance versus the packet error probabilities for optimal and suboptimal solutions for imperfect feedback case, together with the performance of the optimal sequence in the perfect feedback case}
\label{Simu6} 
\end{center}                             
\end{figure}

\subsection{Imperfect Feedback Communication Channel Case} We now consider the case of imperfect packet receipt acknowledgments as described in Section \ref{ternary_subsec} with parameters $\eta=0.4$ and $\delta=0.1$. Let $\lambda$ in (\ref{opt_prob}) be equal to 0.6. The performance of the optimal and suboptimal solutions versus the packet loss probability $\mathbb{P}(\gamma_k=0)=p \in [0.1,0.9]$ is given in Fig. \ref{Simu6}. The performance of the optimal sequence in the case of perfect packet receipt acknowledgments is also shown. We observe that for large packet error probabilities the performance for the suboptimal solution, which is easier to implement, is close to the performance of the optimal solution. 
\vspace*{-0.7cm}
\section{Conclusions and Extensions} \label{conc_sec} This work presents a design methodology for remote estimation of the state of a stable linear stochastic dynamical system, subject to packet dropouts and unreliable acknowledgments. 
The key novelty of this formulation is that the smart sensor decides, at each discrete time instant, whether to transmit either its local state estimate or its local innovation. 
It is shown how to design optimal transmission policies in order to minimize a long term average (infinite-time horizon) cost function as a convex combination of the receiver's expected estimation error variance and the energy needed to transmit the packets. Various computationally efficient suboptimal schemes are presented. For scalar systems,  the optimality of a simple threshold policy in the case of perfect packet receipt acknowledgments is also proved. \\
\indent The analysis of the current paper can be extended to the case of unstable systems with some nontrivial modifications. In order to study unstable systems without feedback control, one can use  the {\em dynamic zoom-in zoom-out quantizer} high rate quantizers as used in \cite{LeongDeyNair_kalmanquant_journal} for decentralized Kalman filtering over bandwidth constrained channels. In case of an unstable system stabilized via feedback control, the approach will likely be different and will possibly use the techniques of linear control design under signal-to-quantization-noise ratio constraints as investigated in \cite{ChiusoLSZ:ecc13,ChiusoLSZ:cdc13}. These and other extensions will be investigated in future work.
\vspace*{-0.4cm}
\bibliographystyle{IEEEtran} 
\renewcommand{\baselinestretch}{0.9}
\bibliography{TCNS-13-0019-Ref}
\vspace*{-0.5cm}
\appendix
\subsection{Proof of Theorem \ref{EEC:thm}}
 For simplicity, denote  $\hat{x}_{k}^s :=  \mathbb{E} [x_{k} |\mathcal Y_{k-1}^s]$. 
 For a given packet loss sequence $\gamma_t,   0 \leq t \leq k-1$, it is evident that $\mathcal Y_{k-1}^r$ is a proper subspace of $\subset Y_{k-1}^s$ due to the fact that both $\hat x^s_{k|k}$ and $\epsilon^s_k$ are linear functions of $y_k$ (see (\ref{KFR:Filter}) and (\ref{InnDef})). Now denote $\hat x_k^r = \E[x_k|\mathcal{Y}_{k-1}^r]$.  Then we have
\begin{align*}
& \E[\hat x^s_k | \mathcal{Y}_{k-1}^r] = \E[\E[x_k| \mathcal Y_{k-1}^s]|\mathcal Y_{k-1}^r] =  \E[x_k|\mathcal{Y}_{k-1}^r] =  \hat x^r_k
\end{align*}
On the other hand,
\begin{align}
&  P_k^{2,2} \equiv \E[(\hat x_k^s- \E[\hat x^s_k | \mathcal{Y}_{k-1}^r])(\hat x_k^s- \E[\hat x^s_k | \mathcal{Y}_{k-1}^r])^T|\mathcal{Y}_{k-1}^r] \notag\\
& \hspace{0.65cm} = \E[(\hat x_k^s-  \hat x^r_k)(\hat x_k^s-  \hat x^r_k)^T|\mathcal{Y}_{k-1}^r]  \notag \\
& \hspace{0.65cm} = \E\big[\big((x_k-\hat x_k^r)-(x_k-\hat x^s_k)\big)\\ & \hspace{0.9cm}\times\big((x_k-\hat x_k^r)-(x_k-\hat x^s_k)\big)^T|\mathcal{Y}_{k-1}^r\big]  \notag \\
&  \hspace{0.65cm} =  P_k^{1,1} + P_s  - 2 \E[(x_k-  \hat x^r_k) (x_k-  \hat x^s_k)^T|\mathcal{Y}_{k-1}^r].
    \label{P22:Eq1}
\end{align}
We note that $\tilde x^s_k:= x_k-  \hat x^s_k$ is orthogonal to $\mathcal Y_{k-1}^s$ and, hence, orthogonal to $\mathcal Y_{k-1}^r$. Therefore, $\E[\hat x_k^s (\tilde x^s_k)^T|\mathcal{Y}_{k-1}^r] =0$ and $E[\hat x_k^r (\tilde x^s_k)^T|\mathcal{Y}_{k-1}^r]=0$ which give
\begin{align*}
& \E[(x_k-  \hat x^r_k) (x_k-  \hat x^s_k)^T|\mathcal{Y}_{k-1}^r]  \\
& \hspace{0.2cm} = \E[\big((x_k- \hat x^s_k)+ (\hat x_k^s- \hat x^r_k)\big) (x_k-  \hat x^s_k)^T|\mathcal{Y}_{k-1}^r] = P_s. 
\end{align*}
This together with (\ref{P22:Eq1}) implies that
$  P_k^{2,2} =  P_k^{1,1} - P_s.
$

In a similar way, we may write 
\begin{align*}
&  P_k^{1,2} 
  =  \E[(x_k-  \hat x^r_k) (\hat x_k^s- \hat x^r_k)^T|\mathcal{Y}_{k-1}^r] \\
&   \hspace{0.65cm} = \E[(x_k-  \hat x^r_k) \big((x_k-\hat x^r_k)-(x_k-\hat x_k^s)\big)^T|\mathcal{Y}_{k-1}^r] \\
& \hspace{0.65cm} =  P_k^{1,1} - P_s.
\hspace{5.7cm} \qed
\end{align*}
\vspace*{-0.5cm}

\subsection{Proof of Lemma \ref{lemmainfstaterec} }
 The total probability
formula\footnote{$ \Prob(A\,|\,B)=\sum_i  \Prob(A,C_i\,|\,B)$} and the chain rule give  
\begin{align}
& \Prob({\bf P}_{k+1},z^k,\nu_k) \! = \! \sum_{\gamma_k} \! \int_{{\bf P}_{k}} \!\! \Prob ({\bf P}_{k+1},{\bf P}_{k},\gamma_k, z^k,\nu_k) d{\bf P}_{k} \notag \\
& \hspace{0.5cm} = \!  \sum_{\gamma_k} \! \int_{{\bf P}_{k}} \!\!\! \Prob ({\bf P}_{k+1}|{\bf P}_{k},\gamma_k, z^k,\nu_k) \Prob ({\bf P}_{k},\gamma_k, z^k,\nu_k)  d{\bf P}_{k} \notag \\
& \hspace{0.5cm} = \! \sum_{\gamma_k} \! \int_{{\bf P}_{k}} \!\!\! \Prob ({\bf P}_{k+1}|{\bf P}_{k},\gamma_k, \nu_k) \Prob ({\bf P}_{k},\gamma_k, z^k,\nu_k)  d{\bf P}_k \label{first_equalit}
\end{align}
where the last equality is because ${\bf P}_{k+1}$ is a function of ${\bf P}_k$, $\gamma_k$ and $\nu_k$ by (\ref{Receiver:RiccatiEq}). However, the chain rule implies that
\begin{align}
& \Prob ({\bf P}_{k},\gamma_k, z^k,\nu_k) \! = \! \Prob ({\bf P}_{k},\gamma_k, z^{k-1},\hat \gamma_k, \nu_{k-1},\nu_k) \notag \\
& = \Prob (\hat \gamma_k | {\bf P}_{k},\gamma_k, z^{k-1},\nu_{k-1},\nu_k) \Prob (\gamma_k | {\bf P}_{k}, z^{k-1}, \nu_{k-1},\nu_k) \notag \\
&  \hspace{0.2cm} \times  \Prob ({\bf P}_k|z^{k-1},\nu_{k-1},\nu_k) \Prob (z^{k-1},\nu_{k-1},\nu_k) \notag \\
& = \Prob (\hat \gamma_k | \gamma_k) \Prob (\gamma_k) \Prob ({\bf P}_k|z^{k-1},\nu_{k-1}) \Prob (z^{k-1},\nu_{k-1},\nu_k).  \label{second_equalit}
\end{align}

Substituting (\ref{second_equalit}) in (\ref{first_equalit}) yields
\begin{align}
& \Prob({\bf P}_{k+1},z^k,\nu_k) = \sum_{\gamma_k} \int_{P_k}  \Big(\Prob ({\bf P}_{k+1}|{\bf P}_{k},\gamma_k, \nu_k)  \Prob (\hat \gamma_k | \gamma_k) \notag \\
& \hspace{0.2cm} \times \Prob (\gamma_k) \Prob ({\bf P}_k|z^{k-1},\nu_{k-1}) \Prob (z^{k-1},\nu_{k-1},\nu_k) \Big) dP_k. \label{simcommon}
\end{align}

On the other hand,
\begin{align}
& \Prob({\bf P}_{k+1}|z^k,\nu_k) = \alpha \times \Prob({\bf P}_{k+1},z^k,\nu_k)  \label{alphaval}
\end{align}
where $\alpha$ is a normalizing constant. Integrating (\ref{alphaval}) with respect to ${\bf P}_{k+1}$ gives 
$\alpha = \big(\int_{{\bf P}_{k+1}} \Prob({\bf P}_{k+1},z^k,\nu_k) d {\bf P}_{k+1}\big)^{-1}.$ However,
\begin{align}
& \int_{{\bf P}_{k+1}} \Prob({\bf P}_{k+1},z^k,\nu_k) d {\bf P}_{k+1} \notag \\
&  =   \int_{{\bf P}_{k+1}} \Big[\sum_{\gamma_k} \int_{{\bf P}_k}  \Big(\Prob ({\bf P}_{k+1}|{\bf P}_k,\gamma_k,\nu_k) \Prob (\hat \gamma_k | \gamma_k)\Prob (\gamma_k) \notag \\
& \hspace{0.2cm} \times \Prob ({\bf P}_k|z^{k-1},\nu_{k-1}) \Prob (z^{k-1},\nu_{k-1},\nu_k) \Big) d{\bf P}_k \Big] d{\bf P}_{k+1}. \label{normval}
\end{align}
By changing the order of integration, we may simplify (\ref{normval}) as
\begin{align}
& \int_{{\bf P}_{k+1}} \Prob({\bf P}_{k+1},z^k,\nu_k) d {\bf P}_{k+1} =  \Prob (z^{k-1},\nu_{k-1},\nu_k) \notag \\
& \hspace{0.3cm}   \times \sum_{\gamma_k} \int_{{\bf P}_k}  \Big( \big( \int_{{\bf P}_{k+1}} \Prob ({\bf P}_{k+1}|{\bf P}_k,\gamma_k,\nu_k) d{\bf P}_{k+1} \big) \Prob (\hat \gamma_k | \gamma_k) \notag \\
&   \hspace{2cm} \times \Prob (\gamma_k)\Prob ({\bf P}_k|z^{k-1},\nu_{k-1})  \Big) d{\bf P}_k \notag \\
& \hspace{0.2cm} =  \Prob (z^{k-1},\nu_{k-1},\nu_k) \sum_{\gamma_k} \Big( \Prob (\hat \gamma_k | \gamma_k) \Prob (\gamma_k)\notag \\
& \hspace{0.3cm}   \times  \big( \int_{{\bf P}_k}\Prob ({\bf P}_k|z^{k-1},\nu_{k-1}) d{\bf P}_k \big) \Big)  \notag \\
& \hspace{0.2cm} =  \Prob (z^{k-1},\nu_{k-1},\nu_k) \sum_{\gamma_k} \Prob (\hat \gamma_k | \gamma_k)\Prob (\gamma_k)
\end{align}
where we used $\int_{{\bf P}_{k+1}} \Prob ({\bf P}_{k+1}|{\bf P}_k,\gamma_k,\nu_k) d{\bf P}_{k+1} =1$ and $\int_{{\bf P}_k}\Prob ({\bf P}_k|z^{k-1},\nu_{k-1}) d{\bf P}_k=1$. Hence, we have
\begin{align}
&\alpha = \Big(\Prob (z^{k-1},\nu_{k-1},\nu_k) \sum_{\gamma_k} \Prob (\hat \gamma_k | \gamma_k)\Prob (\gamma_k)\Big)^{-1}.  \label{alphasim}
\end{align}
Finally, substituting (\ref{simcommon}) and (\ref{alphasim}) in (\ref{alphaval}) gives (\ref{InfState_recur_tot}). \qed
\subsection{Proof of Lemma \ref{Submodul:lem}}
 Let ${\bf P} := \left[
                         \begin{array}{cc}
                           P^{1,1} &  P^{1,1}-P_s \\
                           P^{1,1}-P_s & P^{1,1}-P_s
                         \end{array}
                       \right]$
where $P^{1,1} \geq P_s$ which implies that $ {\bf P}_1 \in \mathbb{S}$. First, note that
\begin{align*}
& F^{1,1}({\bf P},0) = a^2 P^{1,1} + \sigma_w^2 - (1-p) \times \frac{ a^2 K_f^2 P_s^2}{K_f^2P_s+R}
\end{align*}
\begin{align*}
& \textrm{and } F^{1,1}({\bf P},1) = a^2 P^{1,1} + \sigma_w^2 - (1-p) \\
& \hspace{4cm}  \times \frac{ a^2 \big((P^{1,1}-P_s) + K_f P_s\big)^2}{(P^{1,1} - P_s) + K_f^2P_s +R}.
\end{align*}
We denote $g(x) := \frac{\big((x-P_s) + K_f P_s\big)^2}{(x- P_s)+K_f^2P_s+R}$ for $x \geq P_s.$ Let ${\bf P}_1, {\bf P}_2 \in \mathbb{S}$ be such that ${\bf P}_1 \geq {\bf P}_2$, then the inequality (\ref{state-cost_sub}) is equivalent to
\begin{align}
& a^2 (P^{1,1}_1\!\!-\!P^{1,1}_2) \!-\! (1\!-\!p)  a^2 \big(g(P^{1,1}_1) \!-\!g(P^{1,1}_2) \big) \notag \\
& \hspace{5cm} \leq \! a^2 (P^{1,1}_1\!\!-\!P^{1,1}_2). \label{NewInq}
\end{align}

On the other hand,the derivative $g'(x)$ satisfies
\begin{align}
& g'(x) = \frac{2\big((x-P_s) + K_f P_s\big)\big((x- P_s)+K_f^2P_s+R\big)}{\big((x- P_s)+K_f^2P_s+R\big)^2} \notag \\
& \hspace{1.5cm} - \frac{\big((x-P_s) + K_f P_s \big)^2}{\big((x- P_s)+K_f^2P_s+R\big)^2}. \label{der_g}
\end{align}
In the case $x = P_s - K_f P_s$ where either $P_s=0$ or $\sigma_v^2=0$, (\ref{der_g}) yields $g'(x)=0$. Otherwise, dividing the numerator of right hand side of (\ref{der_g}) by the positive value $(x-P_s) + K_f P_s$:
\begin{align*}
& 2 \big( (x- P_s)+K_f^2P_s+R\big) - \big((x-P_s) + K_f P_s \big)  \\
& \hspace{0.5cm} = x-P_s + 2 K_f^2 P_s + 2 R -K_f P_s  \\
& \hspace{0.5cm} = x- P_s + 2K_f^2 (P_s+\sigma_v^2) + 2 \sigma_q^2 - P_s^2/(P_s+\sigma_v^2)
\end{align*}
by the fact that $R = K_f^2 \sigma_v^2 + \sigma_q^2$. Since $K_f = P_s/(P_s+\sigma_v^2)$, 
\begin{align*}
& x- P_s + 2K_f^2 (P_s+\sigma_v^2) + 2 \sigma_q^2 - P_s^2/(P_s+\sigma_v^2) \\
& \hspace{0.5cm}  = x- P_s + P_s^2/(P_s+\sigma_v^2) + 2 \sigma_q^2 \geq 0.
\end{align*}
Therefore, $g'(x)\geq 0$ for $x \geq P_s$. This together with $P^{1,1}_1 \geq P^{1,1}_2 \geq P_s$ implies that the inequality (\ref{NewInq}) is valid. This gives (\ref{state-cost_sub}), thus $F({\bf P}_1, 1) - F({\bf P}_2, 1) \leq F({\bf P}_1, 0) - F({\bf P}_2, 0)$ based on Theorem \ref{EEC:thm} and the fact that for ${\bf P}_1, {\bf P}_2 \in \mathbb{S}$  we have ${\bf P}_1 \geq {\bf P}_2$ if and only if $P_1^{1,1} \geq P_2^{1,1}$. \qed  

\subsection{Proof of Theorem \ref{Submodul:thm}}
Based on the relative value iteration (\ref{RVFI}), define
\begin{align*}
& L_t({\bf P},\nu) = \mathbb{E} \big[\lambda\mathcal L^{1,1} ({\bf P},\gamma,\nu) + (1-\lambda) J(\nu) \big|{\bf P}, \nu\big] \notag \\
& \hspace{0.5cm} + \mathbb{E} \big[V_{t+1} \big(\mathcal L ({\bf P},\gamma,\nu)\big) | {\bf P}, \nu \big] =: L_t^{(1)}({\bf P},\nu) + L_t^{(2)} ({\bf P},\nu) 
\end{align*}
for $0 \leq t \leq T-1$.

{\it Submodularity of $L_t^{(1)}({\bf P},\nu) $}: Lemma \ref{Submodul:lem} implies that $F({\bf P},\nu)= \mathbb{E} \big[\mathcal L({\bf P},\gamma,\nu) \big|{\bf P}, \nu\big] $ and hence $\big(F(P,\nu)\big)^{1,1}= \mathbb{E} \big[\mathcal L^{1,1}({\bf P},\gamma,\nu) \big|{\bf P}, \nu\big]$ are submodular in $({\bf P},\nu)$. It is evident that $\mathbb{E} \big[J(\nu) \big|\nu\big]$ is also submodular in $({\bf P},\nu)$ since it is independent of ${\bf P}$. Therefore, their convex combination $L_t^{(1)}({\bf P},\nu) $ is submodular in $({\bf P},\nu)$. 

{\it Submodularity of $L_t^{(2)}({\bf P},\nu)$}: First we note that both $\mathcal L ({\bf P},\gamma,0) = \mathbb{E} \big[\mathcal L({\bf P},\gamma,\nu) \big|{\bf P}, \nu=0\big] $ and $\mathcal L ({\bf P},\gamma,1)=\mathbb{E} \big[\mathcal L({\bf P},\gamma,\nu) \big|{\bf P}, \nu=1\big] $ given in (\ref{ICFE}) are concave\footnote{The proof of concavity is based on the fact that a function $f(x)$ is concave in $x$ if and only if $f(x_0+th)$ is concave in the scalar $t$ for all $x_0$ and $h$.} and non-decreasing functions in ${\bf P}$ (see Lemma 1 and 2 in \cite{gupta2006stochastic}). This implies that 
$\mathbb{E} \big[\lambda\mathcal L^{1,1} ({\bf P},\gamma,\nu) + (1-\lambda) J(\nu) \big|{\bf P}, \nu=0\big]$, $\mathbb{E} \big[\lambda\mathcal L^{1,1} ({\bf P},\gamma,\nu) + (1-\lambda) J(\nu) \big|{\bf P}, \nu=1\big]$, 
and therefore
\begin{align*}
& \min_{\nu \in \{0,1\}} \Big(\mathbb{E} \big[\lambda\mathcal L^{1,1} ({\bf P},\gamma,\nu) + (1-\lambda) J(\nu) \big|{\bf P}, \nu\big]  \Big)
\end{align*}
are concave and non-decreasing functions of ${\bf P}$ (note that the expectation operator preserves concavity). By induction and the fact that the composition of two non-decreasing concave functions is itself concave and non-decreasing, one can show that the value function $V_t({\bf P})$ in (\ref{perfecfeed_dp:FiniteT}) is a concave and non-decreasing function of ${\bf P}$.

But, the composition of a non-decreasing concave function $V_t(\cdot)$ with a monotonic submodular function $\mathcal L(\cdot,\gamma,\nu)$ is submodular (see part (c) of Proposition 2.3.5 in \cite{simchi2004logic}). Therefore, $
L_t^{(2)}({\bf P},\nu) = \mathbb{E} \big[V_{t+1} \big(\mathcal L ({\bf P},\gamma,\nu)\big) | {\bf P}, \nu \big]$ is submodular in $({\bf P},\nu)$. 

{\it Submodularity of $L_t({\bf P},\nu)$}: The sum of two submodular functions $L_t({\bf P},\nu)=L_t^{(1)}({\bf P},\nu)+L_t^{(2)}({\bf P},\nu)$ is also submodular. 

As a result of submodular functions, for $0 \leq t \leq T-1$, $\arg\!\min_{\nu \in \{0,1\}} L_t({\bf P},\nu)$ is non-decreasing in ${\bf P}$ and hence non-decreasing in $P^{1,1}$. This monotonicity implies the threshold structure (\ref{THRES:Polic}) since the control space has only two elements $\{0,1\}$ (see \cite{topkis2001supermodularity}).\qed.

\end{document}